# Self-Dual Symmetric Polynomials and Conformal Partitions

Leonid G. Fel

School of Physics and Astronomy,
Raymond and Beverly Sackler Faculty of Exact Sciences
Tel Aviv University, Tel Aviv 69978, Israel
e-mail: lfel@ccsg.tau.ac.il

November 3, 2018

**Abstract**

A conformal partition function $\mathcal{P}_n^m(s)$, which arose in the theory of Diophantine equations supplemented with additional restrictions, is concerned with *self-dual symmetric polynomials* – reciprocal $\mathsf{R}_{S_n}^{\{m\}}$ and skew-reciprocal $\mathsf{S}_{S_n}^{\{m\}}$ algebraic polynomials based on the polynomial invariants of the symmetric group $S_n$. These polynomials form an infinite commutative semigroup. Real solutions $\lambda_n(x_i)$ of corresponding algebraic Eqns have many important properties: homogeneity of 1-st order, duality upon the action of the conformal group $\mathsf{W}$, inverting both function $\lambda_n$ and the variables $x_i$, compatibility with trivial solution, *etc*. Making use of the relationship between Gaussian generating function for conformal partitions and Molien generating function for usual restricted partitions we derived the analytic expressions for $\mathcal{P}_n^m(s)$. The unimodality indices for the reciprocal and skew-reciprocal equations were found. The existence of algebraic functions $\lambda_n(x_i)$ invariant upon the action of both the finite group $G \subset S_n$ and conformal group $\mathsf{W}$ is discussed.

# 1 Introduction

The theory of Diophantine equations has a strong relationship with the polynomial invariants of finite groups and different kinds of partitions . The best example is exhibited by the Sylvester theorem [1] about restricted partition number $\mathcal{W}(s, \mathbf{d}^n)$ of positive integer $s$ with respect to the $n$-tuple of positive integers $\mathbf{d}^n = \{d_1, d_2, ..., d_n\}$. $\mathcal{W}(s, \mathbf{d}^n)$ is also a number of sets of positive integer solutions $(x_1, x_2, ..., x_n)$ of the Diophantine equation $\sum_r^n d_r x_r = s$. It is known that $\mathcal{W}(s, \mathbf{d}^n)$ is equal to the coefficient of $t^s$ in the expansion of the generating function

$$M(\mathbf{d}^n, t) = \prod_{r=1}^n \frac{1}{1 - t^{d_r}} = \sum_{s=0}^\infty \mathcal{W}(s, \mathbf{d}^n)\, t^s\ . \tag{1}$$

With the exponents $d_1, d_2, ..., d_n$ becoming the series of integers $1, 2, 3, ..., n$, the number of waves is $n$ and $\mathcal{W}(s, \mathbf{d}^n)$ is usually referred to as a restricted partition number $\mathcal{W}_n(s)$ of $s$ into parts none of which exceeds $n$.

Another definition of $\mathcal{W}(s, \mathbf{d}^n)$ comes from the polynomial invariant of finite reflection groups which are exhausted by irreducible Coxeter groups. Let $M(\mathbf{d}^n, t)$ be a Molien function of such group $G$, $d_r$ are the degrees of basic invariants, and $n$ is the number of basic invariants [2]. Then $\mathcal{W}(s, \mathbf{d}^n)$ gives a number of algebraic independent polynomial invariants of the $s$-degree for group $G$. According to the Sylvester theorem $\mathcal{W}(s, \mathbf{d}^m)$ splits into non-periodic and periodic parts with the periods $\tau\{\mathbf{d}^k\} = \mathsf{lcm}(d_1, d_2, ..., d_k)$. Explicit formulas of $\mathcal{W}(s, \mathbf{d}^n)$ for all irreducible Coxeter groups and a list of $\mathcal{W}_n(s)$ for the first twelve symmetric groups $S_n$ are presented in [3].

There is another kind of partitions which is of high interest. They arise in the theory of Diophantine equations supplemented with additional restrictions, e.g. inequalities. It turns out that such partitions also arise in a special sort of algebraic equations - reciprocal and skew-reciprocal Eqns - based on the polynomial invariants of the symmetric group $S_n$. Recent progress in the self-dual problem of effective isotropic conductivity in two-dimensional (2D) $n$-component regular and random checkerboards [4] has shown the existence of algebraic functions $\lambda_n(x_i)$ which are invariant upon the action of both full permutation group $S_n$, permuting $n$ positive variables $x_i$, and conformal group $\mathsf{W} : \{x_i \to x_i^{-1},\ \lambda_n \to \lambda_n^{-1}\}$. Every coefficient at the term $\lambda_n^{mn-s}$ in the corresponding algebraic equation of the order $mn$ is built out of independent polynomial invariants for symmetric group $S_n$ while their number $\mathcal{P}_n^m(s)$ represents such specific kind of partition function generalizing $\mathcal{W}_n(s)$ in a sense. This function has many remarkable properties. Due to its relation to the conformal group $\mathsf{W}$ we have called $\mathcal{P}_n^m(s)$ a *conformal* partition function.

The paper is organized into 5 Sections. In Section 2 we define the self-dual symmetric polynomials based on the polynomial invariants of the symmetric group $S_n$ and study their algebraic properties. In Section 3 we discuss the existence of other finite groups $G$ distinguished from symmetric $S_n$ that makes it possible to build out the reciprocal $\mathsf{R}_G(\lambda_n, x_i)$ and skew-reciprocal $\mathsf{S}_G(\lambda_n, x_i)$ Eqns on the basis of homogeneous polynomial invariants $I_{d_k}(G)$ of degrees $d_k$. In Section 4 we study the relationship of conformal partition function $\mathcal{P}_n^m(s)$ with a corresponding Diophantine system and consequently with the Gaussian generating function. It permits us to derive the analytic expressions for $\mathcal{P}_n^m(s)$ based on its relation to usual restricted partition function $\mathcal{W}_n(s)$. In Section 5 we derive the compact expression of unimodality indices for the reciprocal and skew-reciprocal polynomials based on the polynomial invariants of the symmetric group $S_n$ as well as of the direct product of symmetric groups.

# 2 Reciprocal and skew-reciprocal symmetric polynomials

The *reciprocal* and *skew-reciprocal* symmetric polynomials are a simple generalization of such polynomials for one variable. Before going forward let us recall some definitions.

Following [5] a polynomial

$$\mathsf{R}_{S_1}^{\{m\}}(\lambda) = \sum_{k=0}^{m} p_k \lambda^k , \quad p_k = p_{m-k} = c_k \qquad (2)$$

will be called a *reciprocal* polynomial. The *skew-reciprocal* polynomial could be defined in a similar manner

$$\mathsf{S}_{S_1}^{\{m\}}(\lambda) = \sum_{k=0}^{m} q_k \lambda^k , \quad q_k = -q_{m-k} = C_k . \qquad (3)$$

The numbers $\mu\left\{\mathsf{R}_{S_1}^{\{m\}}\right\}$ and $\mu\left\{\mathsf{S}_{S_1}^{\{m\}}\right\}$ of independent coefficients $c_k$ and $C_k$ in (2) and (3) are called *unimodality indices* of such polynomials

$$\mu\left\{\mathsf{R}_{S_1}^{\{m\}}\right\} = 1 + \left[\frac{m}{2}\right] , \quad \mu\left\{\mathsf{S}_{S_1}^{\{m\}}\right\} = \left[\frac{m+1}{2}\right] .$$

Now we proceed with the $\mathsf{R}_{S_n}^{\{m\}}$- and $\mathsf{S}_{S_n}^{\{m\}}$- polynomials built upon the basic polynomial invariants $I_{n,r}(x_i)$ of the symmetric group $S_n$

$$I_{n,r} = \sum_{i_1 < i_2 < ... < i_r}^{n} x_{i_1} x_{i_2} ... x_{i_r} , \quad \text{i.e.} \qquad (4)$$

$$I_{n,0} = 1 , \quad I_{n,1} = \sum_{i}^{n} x_i , \quad I_{n,2} = \sum_{i<j}^{n} x_i x_j , \quad ... , \quad I_{n,n-1} = I_{n,n} \sum_{i}^{n} \frac{1}{x_i} , \quad I_{n,n} = \prod_{i}^{n} x_i .$$

$\mathsf{R}_{S_n}^{\{m\}}$- and $\mathsf{S}_{S_n}^{\{m\}}$- polynomials will be defined as a composition of counter-partner monomial terms $T_{\oplus}^{s,l}\left(\begin{smallmatrix}\{m\}\\S_n\end{smallmatrix}\right)$ and $T_{\ominus}^{s,l}\left(\begin{smallmatrix}\{m\}\\S_n\end{smallmatrix}\right)$

$$\mathsf{R}_{S_n}^{[m]}(\lambda_n, x_i) = \sum_{s=0}^{mn} \sum_{l=1}^{\mathcal{P}_n^m(s)} c_{s,l} \left[T_{\oplus}^{s,l}\left(\begin{smallmatrix}\{m\}\\S_n\end{smallmatrix}\right) + T_{\ominus}^{s,l}\left(\begin{smallmatrix}\{m\}\\S_n\end{smallmatrix}\right)\right] ,$$

$$\mathsf{S}_{S_n}^{[m]}(\lambda_n, x_i) = \sum_{s=0}^{mn} \sum_{l=1}^{\mathcal{P}_n^m(s)} C_{s,l} \left[T_{\oplus}^{s,l}\left(\begin{smallmatrix}\{m\}\\S_n\end{smallmatrix}\right) - T_{\ominus}^{s,l}\left(\begin{smallmatrix}\{m\}\\S_n\end{smallmatrix}\right)\right] , \qquad (5)$$

$$T_{\oplus}^{s,l}\left(\begin{smallmatrix}\{m\}\\S_n\end{smallmatrix}\right) = \prod_{r=1}^{n} I_{n,r}^{\alpha_r^{s,l}} \cdot \lambda_n^{mn - \sum_{r=1}^{n} r \alpha_r^{s,l}} , \quad T_{\ominus}^{s,l}\left(\begin{smallmatrix}\{m\}\\S_n\end{smallmatrix}\right) = I_{n,n}^{m - \sum_{r=1}^{n} \alpha_r^{s,l}} \prod_{r=1}^{n} I_{n,n-r}^{\alpha_r^{s,l}} \cdot \lambda_n^{\sum_{r=1}^{n} r \alpha_r^{s,l}} ,$$

where $c_{s,l}$ and $C_{s,l}$ are real numbers, $\left\{\alpha_1^{s,l}, \alpha_2^{s,l}, ..., \alpha_n^{s,l}\right\}$ is a $l$-th tuple of non-negative integers which satisfy

$$\sum_{r=1}^{n} r \alpha_r^{s,l} = s , \quad \sum_{r=1}^{n} \alpha_r^{s,l} \leq m , \qquad (6)$$

and $\mathcal{P}_n^m(s)$ is a number of sets of non-negative integer solutions $\left\{\alpha_1^{s,l}, \alpha_2^{s,l}, ..., \alpha_n^{s,l}\right\}$ of one Diophantine equation and one Diophantine inequality (6), or a number of partitions of $s$ into at most $m$ parts, each $\leq n$. The definition (5) preserves *a constitutive relation*

$$T_\oplus^{s,l}\left(\begin{array}{c}\{m\}\\S_n\end{array}\right) \times T_\ominus^{s,l}\left(\begin{array}{c}\{m\}\\S_n\end{array}\right) = I_{n,n}^{m-\sum_{r=1}^n \alpha_r^{s,l}} \prod_{r=1}^n (I_{n,r} I_{n,n-r})^{\alpha_r^{s,l}} \cdot \lambda_n^{mn} \tag{7}$$

It is worth representing the $\mathsf{R}_{S_n}^{\{m\}}$- and $\mathsf{S}_{S_n}^{\{m\}}$- polynomials as a four-part decomposition

$$\mathsf{R}_{S_n}^{\{m\}}(\lambda_n, x_i) = \mathsf{T1}_{S_n}^{\{m\}}(\lambda_n, x_i) + \mathsf{T2}_{S_n}^{\{m\}}(\lambda_n, x_i) + \mathsf{T3}_{S_n}^{\{m\}}(\lambda_n, x_i) + \mathsf{T4}_{S_n}^{\{m\}}(\lambda_n, x_i), \tag{8}$$

$$\mathsf{S}_{S_n}^{\{m\}}(\lambda_n, x_i) = \mathsf{T1}_{S_n}^{\{m\}}(\lambda_n, x_i) + \mathsf{T2}_{S_n}^{\{m\}}(\lambda_n, x_i) - \mathsf{T3}_{S_n}^{\{m\}}(\lambda_n, x_i) - \mathsf{T4}_{S_n}^{\{m\}}(\lambda_n, x_i), \tag{9}$$

$$\mathsf{T1}_{S_n}^{\{m\}}(\lambda_n, x_i) = \sum_{s=0}^{\left[\frac{mn-1}{2}\right]} \lambda_n^{mn-s} \sum_{l=1}^{\mathcal{P}_n^m(s)} c_{s,l} \prod_{r=1}^n I_{n,r}^{\alpha_r^{s,l}} = c_{0,1} \lambda_n^{mn} + c_{1,1} I_{n,1} \cdot \lambda_n^{mn-1} + \tag{10}$$

$$(c_{2,1} I_{n,1}^2 + c_{2,2} I_{n,2}) \cdot \lambda_n^{mn-2} + (c_{3,1} I_{n,1}^3 + c_{3,2} I_{n,1} I_{n,2} + c_{3,3} I_{n,3}) \cdot \lambda_n^{mn-3} +$$

$$(c_{4,1} I_{n,1}^4 + c_{4,2} I_{n,1}^2 I_{n,2} + c_{4,3} I_{n,2}^2 + c_{4,4} I_{n,1} I_{n,3} + c_{4,5} I_{n,4}) \cdot \lambda_n^{mn-4} + \ldots,$$

$$\mathsf{T2}_{S_n}^{\{m\}}(\lambda_n, x_i) = \lambda_n^{\frac{mn}{2}} \sum_{l=1}^{\mathcal{P}_n^m(\frac{mn}{2})} c_{s^*,l} \prod_{r=1}^n I_{n,r}^{\alpha_r^{s^*,l}} \cdot \theta(mn), \quad s^* = \frac{mn}{2},$$

$$\mathsf{T3}_{S_n}^{\{m\}}(\lambda_n, x_i) = \lambda_n^{\frac{mn}{2}} \sum_{l=1}^{\mathcal{P}_n^m(\frac{mn}{2})} c_{s^*,l} I_{n,n}^{(m-\sum_{r=1}^n \alpha_r^{s^*,l})} \prod_{r=1}^n I_{n,n-r}^{\alpha_r^{s^*,l}} \cdot \theta(mn),$$

$$\mathsf{T4}_{S_n}^{\{m\}}(\lambda_n, x_i) = \sum_{s=0}^{\left[\frac{mn-1}{2}\right]} \lambda_n^s \sum_{l=1}^{\mathcal{P}_n^m(s)} c_{s,l} I_{n,n}^{(m-\sum_{r=1}^n \alpha_r^{s,l})} \prod_{r=1}^n I_{n,n-r}^{\alpha_r^{s,l}} = c_{0,1} I_{n,n}^m + c_{1,1} I_{n,n-1} I_{n,n}^{m-1} \cdot \lambda_n +$$

$$(c_{2,1} I_{n,n-1}^2 + c_{2,2} I_{n,n-2} I_{n,n}) I_{n,n}^{m-2} \cdot \lambda_n^2 + (c_{3,1} I_{n,n-1}^3 + c_{3,2} I_{n,n-2} I_{n,n-1} I_{n,n} +$$

$$c_{3,3} I_{n,n-3} I_{n,n}^2) I_{n,n}^{m-3} \cdot \lambda_n^3 + (c_{4,1} I_{n,n-1}^4 + c_{4,2} I_{n,n-1}^2 I_{n,n-2} I_{n,n} + c_{4,3} I_{n,n-2}^2 I_{n,n}^2 +$$

$$c_{4,4} I_{n,n-1} I_{n,n-3} I_{n,n}^2 + c_{4,5} I_{n,n-4} I_{n,n}^3) I_{n,n}^{m-4} \cdot \lambda_n^4 + \ldots.$$

where $\theta(mn)$ denotes a step-function : $\theta(2k) = 1$, $\theta(2k+1) = 0$.

We also define the corresponding unimodality indices $\mu\left\{\mathsf{R}_{S_n}^{\{m\}}\right\}$ and $\mu\left\{\mathsf{S}_{S_n}^{\{m\}}\right\}$ as numbers of independent coefficients $c_{k,l}$ in (8) and (9)

$$\mu\left\{\mathsf{R}_{S_n}^{\{m\}}\right\} = \mu\left\{\mathsf{T1}_{S_n}^{\{m\}}\right\} + \mu\left\{\mathsf{T2}_{S_n}^{\{m\}} + \mathsf{T3}_{S_n}^{\{m\}}\right\}, \quad \mu\left\{\mathsf{T2}_{S_n}^{\{m\}}\right\} = \mu\left\{\mathsf{T3}_{S_n}^{\{m\}}\right\},$$

$$\mu\left\{\mathsf{S}_{S_n}^{\{m\}}\right\} = \mu\left\{\mathsf{T1}_{S_n}^{\{m\}}\right\} + \mu\left\{\mathsf{T2}_{S_n}^{\{m\}} - \mathsf{T3}_{S_n}^{\{m\}}\right\}, \quad \mu\left\{\mathsf{T1}_{S_n}^{\{m\}}\right\} = \mu\left\{\mathsf{T4}_{S_n}^{\{m\}}\right\},$$

$$\mu\left\{\mathsf{T2}_{S_n}^{\{m\}} + \mathsf{T3}_{S_n}^{\{m\}}\right\} > \mu\left\{\mathsf{T2}_{S_n}^{\{m\}} - \mathsf{T3}_{S_n}^{\{m\}}\right\}. \tag{11}$$

The inequality (11) reflects the fact that some terms in (9) annihilate with their counterpartners. The following is an instructive example

$$\mathsf{R}_{S_4}^{\{2\}} = c_{0,1} \lambda_4^8 + c_{1,1} I_{4,1} \lambda_4^7 + (c_{2,1} I_{4,1}^2 + c_{2,2} I_{4,2}) \lambda_4^6 + (c_{3,1} I_{4,1} I_{4,2} + c_{3,2} I_{4,3}) \lambda_4^5 +$$

$$(c_{4,1} I_{4,1} I_{4,3} + c_{4,2} I_{4,2}^2 + c_{4,3} I_{4,4}) \lambda_4^4 +$$

$$(c_{3,1} I_{4,2} I_{4,3} + c_{3,2} I_{4,1} I_{4,4}) \lambda_4^3 + (c_{2,1} I_{4,3}^2 + c_{2,2} I_{4,2} I_{4,4}) \lambda_4^2 + c_{1,1} I_{4,3} I_{4,4} \lambda_4 + c_{0,1} I_{4,4}^2,$$

$$\mathsf{S}_{S_4}^{\{2\}} = C_{0,1} \lambda_4^8 + C_{1,1} I_{4,1} \lambda_4^7 + (C_{2,1} I_{4,1}^2 + C_{2,2} I_{4,2}) \lambda_4^6 + (C_{3,1} I_{4,1} I_{4,2} + C_{3,2} I_{4,3}) \lambda_4^5 -$$

$$(C_{3,1} I_{4,2} I_{4,3} + C_{3,2} I_{4,1} I_{4,4}) \lambda_4^3 - (C_{2,1} I_{4,3}^2 + C_{2,2} I_{4,2} I_{4,4}) \lambda_4^2 - C_{1,1} I_{4,3} I_{4,4} \lambda_4 - C_{0,1} I_{4,4}^2,$$

that gives

$$\mu\left\{\mathsf{T1}_{S_4}^{\{2\}}\right\} = 6,\ \mu\left\{\mathsf{T2}_{S_4}^{\{2\}} + \mathsf{T3}_{S_4}^{\{2\}}\right\} = 3,\ \mu\left\{\mathsf{T2}_{S_4}^{\{2\}} - \mathsf{T3}_{S_4}^{\{2\}}\right\} = 0,\ \mu\left\{\mathsf{R}_{S_4}^{\{2\}}\right\} = 9,\ \mu\left\{\mathsf{S}_{S_4}^{\{2\}}\right\} = 6\ .$$

It is often useful to consider *monic* polynomials which have $c_{0,1} = 1$, or $C_{0,1} = 1$. Their unimodality indices are $\mu\left\{\mathsf{R}_{S_n}^{\{m\}}\right\} - 1$ and $\mu\left\{\mathsf{S}_{S_n}^{\{m\}}\right\} - 1$.

We consider some important properties of the $\mathsf{R}_{S_n}^{\{m\}}$- and $\mathsf{S}_{S_n}^{\{m\}}$- polynomials.

**Proposition 1**
$\mathsf{R}_{S_n}^{\{m\}}$- and $\mathsf{S}_{S_n}^{\{m\}}$- polynomials form an infinite commutative semigroup :

$$\mathsf{R}_{S_n}^{\{m_1\}} \times \mathsf{R}_{S_n}^{\{m_2\}} = \mathsf{S}_{S_n}^{\{m_1\}} \times \mathsf{S}_{S_n}^{\{m_2\}} = \mathsf{R}_{S_n}^{\{m_1+m_2\}};\ \mathsf{R}_{S_n}^{\{m_1\}} \times \mathsf{S}_{S_n}^{\{m_2\}} = \mathsf{S}_{S_n}^{\{m_1\}} \times \mathsf{R}_{S_n}^{\{m_2\}} = \mathsf{S}_{S_n}^{\{m_1+m_2\}} \quad (12)$$

$\mathsf{R}_{S_n}^{\{m\}}$ - polynomials form a sub-semigroup of index two in the above semigroup.

*Proof.* In order to prove (12) let us consider two skew-reciprocal polynomials $\mathsf{S}_{S_n}^{\{m_1\}}$ and $\mathsf{S}_{S_n}^{\{m_2\}}$. Multiplying them we will get two kinds of generic terms $T_{3\oplus}^{s,l}$, $T_{4\oplus}^{s,l}$ with their counter-partners $T_{3\ominus}^{s,l}$, $T_{4\ominus}^{s,l}$ :

$$T_{3\oplus}^{s,l} = T_{1\oplus}^{s,l}\left(\mathsf{S}_{S_n}^{\{m_1\}}\right) \times T_{2\oplus}^{s,l}\left(\mathsf{S}_{S_n}^{\{m_2\}}\right) = \prod_{r=1}^{n} I_{n,r}^{\alpha_r^{s,l}+\beta_r^{s,l}} \cdot \lambda_n^{(m_1+m_2)n - \sum_{r=1}^{n} r\,(\alpha_r^{s,l}+\beta_r^{s,l})}\ ,$$

$$T_{3\ominus}^{s,l} = T_{1\ominus}^{s,l}\left(\mathsf{S}_{S_n}^{\{m_1\}}\right) \times T_{2\ominus}^{s,l}\left(\mathsf{S}_{S_n}^{\{m_2\}}\right) = I_{n,n}^{(m_1+m_2) - \sum_{r=1}^{n}(\alpha_r^{s,l}+\beta_r^{s,l})} \prod_{r=1}^{n} I_{n,n-r}^{\alpha_r^{s,l}+\beta_r^{s,l}} \cdot \lambda_n^{\sum_{r=1}^{n} r\,(\alpha_r^{s,l}+\beta_r^{s,l})}\ ,$$

$$T_{4\oplus}^{s,l} = T_{1\oplus}^{s,l}\left(\mathsf{S}_{S_n}^{\{m_1\}}\right) \times T_{2\ominus}^{s,l}\left(\mathsf{S}_{S_n}^{\{m_2\}}\right) = I_{n,n}^{m_2 - \sum_{r=1}^{n} \beta_r^{s,l}} \prod_{r=1}^{n} I_{n,r}^{\alpha_r^{s,l}} I_{n,n-r}^{\beta_r^{s,l}} \cdot \lambda_n^{m_1 n + \sum_{r=1}^{n} r\,(\beta_r^{s,l} - \alpha_r^{s,l})}\ ,$$

$$T_{4\ominus}^{s,l} = T_{1\ominus}^{s,l}\left(\mathsf{S}_{S_n}^{\{m_1\}}\right) \times T_{2\oplus}^{s,l}\left(\mathsf{S}_{S_n}^{\{m_2\}}\right) = I_{n,n}^{m_1 - \sum_{r=1}^{n} \alpha_r^{s,l}} \prod_{r=1}^{n} I_{n,r}^{\beta_r^{s,l}} I_{n,n-r}^{\alpha_r^{s,l}} \cdot \lambda_n^{m_2 n + \sum_{r=1}^{n} r\,(\alpha_r^{s,l} - \beta_r^{s,l})}\ .$$

We obviously arrive at the constitutive relation claimed in (7) for polynomial $\mathsf{R}_{S_n}^{\{m_1+m_2\}}$

$$T_{3\oplus}^{s,l} \times T_{3\ominus}^{s,l} = T_{4\oplus}^{s,l} \times T_{4\ominus}^{s,l} = I_{n,n}^{(m_1+m_2) - \sum_{r=1}^{n}(\alpha_r^{s,l}+\beta_r^{s,l})} \prod_{r=1}^{n}(I_{n,r} I_{n,n-r})^{\alpha_r^{s,l}+\beta_r^{s,l}} \cdot \lambda_n^{(m_1+m_2)n}$$

and the number $\mathcal{P}_n^{m_1+m_2}(s)$ of the terms composed with equal power $\lambda^s$ is determined by Diophantine systems

$$\alpha_r^{s,l}\,,\ \beta_r^{s,l} \geq 0\,,\quad \sum_{r=1}^{n} r\,(\alpha_r^{s,l} + \beta_r^{s,l}) = s\,,\quad \sum_{r=1}^{n}(\alpha_r^{s,l} + \beta_r^{s,l}) \leq m_1 + m_2\ .$$

Thus, the monomials $T_{3\oplus}^{s,l}$ and $T_{3\ominus}^{s,l}$ are the counter-partners as well as the monomials $T_{4\oplus}^{s,l}$ and $T_{4\ominus}^{s,l}$. They compose the terms of $\mathsf{R}_{S_n}^{\{m_1+m_2\}}$- polynomials. The remaining multiplication rules (12) could be proved in the similar way. The commutativity follows immediately. The second part of Proposition follows from the multiplication rules (12). ∎

## 2.1 Homogeneity, duality and compatibility

The corresponding algebraic equations $\mathsf{R}_{S_n}^{\{m\}}(\lambda_n, x_i) = 0$ and $\mathsf{S}_{S_n}^{\{m\}}(\lambda_n, x_i) = 0$ will be called $\mathsf{R}_{S_n}^{\{m\}}$- and $\mathsf{S}_{S_n}^{\{m\}}$- equations respectively. Let $\lambda_n(x_i)$ be a real solution of one of the $\mathsf{R}_{S_n}^{\{m\}}$- or $\mathsf{S}_{S_n}^{\{m\}}$- Eqns. Then we are able to prove the following

**Proposition 2**

$\lambda_n(x_i)$ is a homogeneous of 1-st order algebraic function, which is invariant upon the action of the symmetric group $S_n$, permuting n non-negative variables $x_i$

$$\lambda_n(k \cdot x_i) = k \cdot \lambda_n(x_i), \quad \lambda_n(\widehat{\mathcal{D}}\{x_1, x_2, ..., x_n\}) = \lambda_n(x_1, x_2, ..., x_n), \tag{13}$$

where $\widehat{\mathcal{D}}$ is a permutation operator of the indices $\{1, 2, ..., n\}$.

*Proof.* The $S_n$-permutation invariance of $\lambda_n(x_i)$ does exist due to (4), (5). The homogeneity of 1-st order could be established by calculation of the degree of the generic terms

$$\deg\left[T_{\oplus}^{s,l}\binom{\{m\}}{S_n}(\lambda_n, x_i)\right] = \deg\left[\lambda_n^{mn - \sum_{r=1}^n r\, \alpha_r^{s,l}}\right] + \deg\left[\prod_{r=1}^n I_{n,r}^{\alpha_r^{s,l}}\right] = mn,$$

$$\deg\left[T_{\ominus}^{s,l}\binom{\{m\}}{S_n}(\lambda_n, x_i)\right] = \deg\left[\lambda_n^{\sum_{r=1}^n r\, \alpha_r^{s,l}}\right] + \deg\left[I_{n,n}^{m - \sum_{r=1}^n \alpha_r^{s,l}}\right] + \deg\left[\prod_{r=1}^n I_{n,n-r}^{\alpha_r^{s,l}}\right] = mn,$$

which is independent of index $l$. ∎

**Proposition 3**

$\lambda_n(x_i)$ is invariant upon the action of the conformal group $\mathsf{W}$, inverting both a function $\lambda_n$ and the variables $x_i$

$$\lambda_n(x_1, x_2, ..., x_n) \times \lambda_n\left(\frac{1}{x_1}, \frac{1}{x_2}, ..., \frac{1}{x_n}\right) = 1. \tag{14}$$

*Proof.* This property follows when we note that (14) is equivalent to e.g. $\mathsf{S}_{S_n}^{\{m\}}$- Eqn for $\frac{1}{\lambda_n}(x_1^{-1}, x_2^{-1}, ..., x_n^{-1})$. Let us apply a conformal transformation

$$\widehat{\mathsf{W}} = \left\{x_i \leftrightarrow \frac{1}{x_i},\ \lambda_n \leftrightarrow \frac{1}{\lambda_n}\right\} \quad \Rightarrow \quad \mathsf{W}: \left\{I_{n,r} \longleftrightarrow \frac{I_{n,n-r}}{I_{n,n}}\right\} \tag{15}$$

to the counter-partner monomial terms $T_{\oplus}^{s,l}\binom{\{m\}}{S_n}$ and $T_{\ominus}^{s,l}\binom{\{m\}}{S_n}$

$$T_{\oplus}^{s,l}\binom{\{m\}}{S_n} \xrightarrow{\widehat{W}} U_{\oplus}^{s,l}\binom{\{m\}}{S_n} = I_{n,n}^{-\sum_{r=1}^n \alpha_r^{s,l}} \prod_{r=1}^n I_{n,n-r}^{\alpha_r^{s,l}} \cdot \lambda_n^{\sum_{r=1}^n r\, \alpha_r^{s,l} - mn}$$

$$T_{\ominus}^{s,l}\binom{\{m\}}{S_n} \xrightarrow{\widehat{W}} U_{\ominus}^{s,l}\binom{\{m\}}{S_n} = I_{n,n}^{-m} \prod_{r=1}^n I_{n,r}^{\alpha_r^{s,l}} \cdot \lambda_n^{-\sum_{r=1}^n r\, \alpha_r^{s,l}},$$

The new terms $U_{\oplus}^{s,l}\binom{\{m\}}{S_n}$ and $U_{\ominus}^{s,l}\binom{\{m\}}{S_n}$ preserve the constitutive relation (7) up to the inessential multiplier

$$(\lambda_n I_{n,n})^{2mn} \times U_{\oplus}^{s,l}\binom{\{m\}}{S_n} \times U_{\ominus}^{s,l}\binom{\{m\}}{S_n} = I_{n,n}^{m - \sum_{r=1}^n \alpha_r^{s,l}} \prod_{r=1}^n (I_{n,r} I_{n,n-r})^{\alpha_r^{s,l}} \cdot \lambda_n^{mn}.$$

It is easy to make sure that the solution $\lambda_n(x_i)$ of both $\mathsf{R}_{S_n}^{\{m\}}(\lambda_n, x_i)$ and $\mathsf{S}_{S_n}^{\{m\}}(\lambda_n, x_i)$ Eqns are conformally invariant. Indeed, due to the transformation (15) for invariants we have

$$\mathsf{R}_{S_n}^{\{m\}}(\lambda_n, x_i) \xrightarrow{\widehat{W}} \mathsf{R}_{S_n}^{\{m\}}\left(\frac{1}{\lambda_n}, \frac{1}{x_i}\right) = (\lambda_n I_{n,n})^{-m} \cdot \mathsf{R}_{S_n}^{\{m\}}(\lambda_n, x_i) , \qquad (16)$$

$$\mathsf{S}_{S_n}^{\{m\}}(\lambda_n, x_i) \xrightarrow{\widehat{W}} \mathsf{S}_{S_n}^{\{m\}}\left(\frac{1}{\lambda_n}, \frac{1}{x_i}\right) = -(\lambda_n I_{n,n})^{-m} \cdot \mathsf{S}_{S_n}^{\{m\}}(\lambda_n, x_i) .$$

The above given Eqns prove the existence of algebraic functions $\lambda_n(x_i)$ which are invariant upon the action both symmetric group $S_n$, permuting $n$ positive variables $x_i$, and conformal group $\mathsf{W}$, inverting both function $\lambda_n$ and the variables $x_i$. ∎

It is quite interesting that these two properties (13) and (14) imposed on the algebraic function $\lambda_n(x_i)$ are also sufficient to specify a certain sort of algebraic Eqns which give rise to such function.

**Proposition 4**
*Let $\lambda_n(x_i)$ be a homogeneous of 1-st order algebraic function, which is invariant upon the action of both symmetric group $S_n$ (13) and conformal group $\mathsf{W}$ (14). Then necessarily $\lambda_n(x_i)$ is a real solution of only $\mathsf{R}_{S_n}^{\{m\}}(\lambda_n, x_i)$ - or $\mathsf{S}_{S_n}^{\{m\}}(\lambda_n, x_i)$ Eqns.*

*Proof.* Let us consider a generic form of algebraic Eqn which satisfies (13), (14). Since the basic invariant of the highest order is $I_{n,n}$, the degree of Eqn $\mathcal{F}(\lambda_n, x_i) = 0$, which preserves the homogeneity, has $n$ as a divisor

$$\mathcal{F}(\lambda_n, x_i) = \sum_{s=0}^{mn} A_s(I_{n,r}) \cdot \lambda_n^{mn-s} , \quad A_s(I_{n,r}) = \sum_{l=1}^{\mathcal{W}_n(s)} c_{s,l} \cdot \prod_{r=1}^{n} I_{n,r}^{\alpha_r^{s,l}} , \quad \sum_{r=1}^{n} r\, \alpha_r^{s,l} = s , \qquad (17)$$

where $mn$ is an order of the Eqn, $\mathcal{W}_n(s)$ is a number of algebraic independent polynomial invariants of $s$-degree for group $S_n$ [3] and the variable $l$ numerates different tuples of positive integers solutions $\left\{\alpha_1^{s,l}, \alpha_2^{s,l}, ..., \alpha_n^{s,l}\right\}$ of the Diophantine equation in (17). Applying conformal transformation (15) to Eqn (17) and multiplying by $I_{n,n}^m \lambda_n^{mn}$ we get a new Eqn $\widetilde{\mathcal{F}}(\lambda_n, x_i) = 0$

$$\widetilde{\mathcal{F}}(\lambda_n, x_i) = \sum_{s=0}^{mn} B_s(I_{n,r}) \cdot \lambda_n^s , \quad B_s(I_{n,r}) = \sum_{l=1}^{\mathcal{W}_n(s)} c_{s,l} \cdot I_{n,n}^{(m-\sum_{r=1}^{n}\alpha_r^{s,l})} \cdot \prod_{r=1}^{n} I_{n,n-r}^{\alpha_r^{s,l}} , \quad \sum_{r=1}^{n} \alpha_r^{s,l} \leq m .$$

To satisfy the $\mathcal{F} \equiv \pm\widetilde{\mathcal{F}}$ coincidence requirement the corresponding terms should be equated. It diminishes the number of permitted tuples $\left\{\alpha_1^{s,l}, \alpha_2^{s,l}, ..., \alpha_n^{s,l}\right\}$ up to $\mathcal{P}_n^m(s)$ in accordance with Diophantine system (6). Rearranging the summation and product in (17)

$$\mathcal{F}(\lambda_n, x_i) = \sum_{s'=0}^{mn} A_{s'}(I_{n,r}) \cdot \lambda_n^s , \quad A_{s'}(I_{n,r}) = \sum_{l=1}^{\mathcal{P}_n^m(s')} c_{s',l} \cdot \prod_{r=0}^{n-1} I_{n,n-r}^{\alpha_{n-r}^{s',l}} , \quad s' = mn - s$$

and comparing it with $\widetilde{\mathcal{F}}(\lambda_n, x_i)$ we will get $\sum_{s=0}^{mn} \mathcal{P}_n^m(s)$ identities

$$\sum_{l=1}^{\mathcal{P}_n^m(s)} c_{s,l} \cdot I_{n,n}^{(m-\sum_{r=1}^{n}\alpha_r^{s,l})} \cdot \prod_{r=1}^{n} I_{n,n-r}^{\alpha_r^{s,l}} = \pm \sum_{l=1}^{\mathcal{P}_n^m(s')} c_{s',l} \cdot \prod_{r=0}^{n-1} I_{n,n-r}^{\alpha_{n-r}^{s',l}} , \quad \mathcal{P}_n^m(s') = \mathcal{P}_n^m(s) ,$$

where the latter relation for conformal partition $\mathcal{P}_n^m(s')$ follows from the Cayley-Sylvester theorem (see Section 4). Thus, due to algebraic independence of basic invariants we arrive at the following relations

$$\alpha_r^{s,l} = \alpha_{n-r}^{s',l} \,, \quad m - \sum_{r=1}^{n} \alpha_r^{s,l} = \alpha_n^{s',l} \,, \quad c_{s,l} = \pm c_{s',l} \,. \tag{18}$$

Therefore one can choose the counter-partners $T_\oplus^{s,l}\binom{\{m\}}{S_n}(\lambda_n, x_i)$ and $T_\ominus^{s,l}\binom{\{m\}}{S_n}(\lambda_n, x_i)$ as was done in (5). ■

**Proposition 5**
The real solution $\lambda_n(x_i)$ of $\mathsf{S}_{S_n}^{\{m\}}(\lambda_n, x_i)$ Eqn is compatible with equating all variables $x_i$

$$\lambda_n(x, x, ..., x) = x \,. \tag{19}$$

*Proof.* Let us calculate the difference between two counter-partner monomial terms $T_\oplus^{s,l}\binom{\{m\}}{S_n}(\lambda_n, x_i)$ and $T_\ominus^{s,l}\binom{\{m\}}{S_n}(\lambda_n, x_i)$, defined in (5), when $\lambda_n = x_i = x$:

$$0 = \mathsf{S}_{S_n}^{\{m\}}(x,x) = \sum_{s=0}^{mn} \sum_{l=1}^{\mathcal{P}_n^m(s)} C_{s,l} \left[ T_\oplus^{s,l}\binom{\{m\}}{S_n}(x,x) - T_\ominus^{s,l}\binom{\{m\}}{S_n}(x,x) \right] =$$

$$\sum_{s=0}^{mn} \sum_{l=1}^{\mathcal{P}_n^m(s)} C_{s,l} \left[ x^{nm - \sum_{r=1}^{n} r\alpha_r^{s,l}} \cdot x^{\sum_{r=1}^{n} r\alpha_r^{s,l}} - x^{nm - n\sum_{r=1}^{n} \alpha_r^{s,l}} \cdot x^{\sum_{r=1}^{n} r\alpha_r^{s,l} + \sum_{r=1}^{n} (n-r)\alpha_r^{s,l}} \right] \prod_{r=1}^{[\frac{n}{2}]} \binom{n}{r}^{\alpha_r^{s,l}}$$

that proves the Proposition. ■

## 2.2 Skew-reciprocal Eqns: exact solution and universal bounds

In this Section we will discuss some kinds of skew-reciprocal Eqns which allow to find their exact solutions $\lambda_n(x_i)$ or at least establish the upper and lower bounds for $\lambda_n(x_i)$. It is easy to see from (3) that monic $\mathsf{S}_{S_1}^{\{1\}}$- polynomial is the exact divisor of any $\mathsf{S}_{S_1}^{\{m\}}$- polynomial. What is less trivial that the similar property also have $\mathsf{S}_{S_2}^{\{m\}}$-polynomials.

**Proposition 6**

$$\mathsf{S}_{S_2}^{\{m\}} = 0 \mod \left( \mathsf{S}_{S_2}^{\{1\}} \right) \,. \tag{20}$$

*Proof.* According to (9) the monic $\mathsf{S}_{S_2}^{\{1\}}$- Eqn is the following: $\mathsf{S}_{S_2}^{\{1\}} = \lambda_2^2 - I_{2,2}$, while the monic $\mathsf{S}_{S_2}^{\{m\}}$- Eqn could be represented as

$$\mathsf{S}_{S_2}^{\{m\}} = \lambda_2^{2m} - I_{2,2}^m + C_{1,1} I_{2,1} \lambda_2 \left( \lambda_2^{2m-2} - I_{2,2}^{m-1} \right) + (C_{2,1} I_{2,1}^2 + C_{2,2} I_{2,2}) \lambda_2^2 \left( \lambda_2^{2m-4} - I_{2,2}^{m-2} \right) + ...$$

where every term in the above expansion has factor $\lambda_2^{2m-2k} - I_{2,2}^{m-k}$, $0 \leq k \leq m - 1$. Thus, $\lambda_2(x_1, x_2) = \sqrt{x_1 x_2}$ is the common solution of $\mathsf{S}_{S_2}^{\{m\}}$- Eqn. This proves (20). ■

It turns out that $n = 2$ is also the last exactly solvable case when the unique function $\lambda_n(x_i)$ is the common solution of all $\mathsf{S}_{S_n}^{\{m\}}$- Eqns. This property disappears for $n \geq 3$ that reflects the growth of unimodality index, i.e. $\mu\left\{\mathsf{S}_{S_1}^{\{1\}}\right\} = \mu\left\{\mathsf{S}_{S_2}^{\{1\}}\right\} = 1$ and $\mu\left\{\mathsf{S}_{S_3}^{\{1\}}\right\} = 2$. Nevertheless we can state and prove that for certain non-trivially related $x_i$ this property survives.

**Proposition 7**

If the polynomial invariants $I_{n,1}, ..., I_{n,n}$ satisfy $r$ identities simultaneously

$$\sqrt[n]{I_{n,n}} = \sqrt[n-2r]{\frac{I_{n,n-r}}{I_{n,r}}}, \quad 2r \leq n, \tag{21}$$

then $\lambda_n(x_i) = \sqrt[n]{I_{n,n}}$ is the common root of all $\mathsf{S}_{S_n}^{\{m\}}$- Eqns.

*Proof.* Equating two counter-partners in (5) and putting there $\lambda_n^n = I_{n,n}$ we arrive at Eqn

$$\sum_{r=1}^{n} \alpha_r^{s,l} \left[ \left(1 - \frac{2r}{n}\right) \ln I_{n,n} + \ln I_{n,r} - \ln I_{n,n-r} \right] = 0,$$

which is independent of $m$. For $n = 1, 2$ the relation (21) is satisfied identically: $n = 1, I_1 = I_1$; $n = 2, I_{2,1}^2 = I_{2,1}^2$. ∎

We will show now that (21) is equivalent to the following statement:

**Proposition 8**

If there exists such $\lambda_*$, that for every $x_i$ one can point out such $x_j$, which satisfies

$$x_i \cdot x_j = \lambda_{n*}^2, \tag{22}$$

then $\lambda_n(x_i) = \sqrt[n]{I_{n,n}}$ is the common root of all $\mathsf{S}_{S_n}^{\{m\}}$- Eqns.

*Proof.* (22)⇒(21). This property follows when we note that (22) yields $I_{n,n} = \lambda_{n*}^n$ and the relation (21) is transformed into $I_{n,r} = I_{n,n-r}\lambda_{n*}^{2r-n}$ or

$$\sum_{i_1<i_2<...<i_r}^{n} x_{i_1}x_{i_2}...x_{i_r} = \sum_{i_1<i_2<...<i_r}^{n} \frac{\lambda_{n*}^{2r}}{x_{i_1}x_{i_2}...x_{i_r}} = \sum_{j_1<j_2<...<j_r}^{n} x_{j_1}x_{j_2}...x_{j_r}. \tag{23}$$

Since the result of summation in (23) does not depend on the summation order, we arrive at the conclusion that (23) is the identity.

(22)⇐(21). Without assumption (22) we can rewrite (21) as follows

$$J_{n,k} = I_{n,k}, \quad 1 \leq k \leq r \leq \left[\frac{n}{2}\right]; \quad J_{n,r} = \sum_{i_1<i_2<...<i_r}^{n} y_{i_1}y_{i_2}...y_{i_r}, \quad y_i = \frac{\sqrt[n]{I_{n,n}^2}}{x_i}, \tag{24}$$

where $y_i$ compose another set $\{y_i\}$ of $n$ non-negative real variables, which generate $n$ basic invariants $J_{n,r}$ in a similar manner (4). The relationship between $\{y_i\}$ and $\{x_i\}$ is defined by $r$ equations (24). We will show that in fact we have $n$ non-linear Eqns, which completely determine the $\{y_i\} = \{x_i\}$ coincidence up to the permutation $\widehat{\mathcal{D}}$ of the variables $y_i$ and $x_i$ in each set separately. Indeed, let us determine a following invariant $J_{n,r+1}$, which is absent in (24)

$$J_{n,r+1} = \sum_{i_1<i_2<...<i_{r+1}}^{n} y_{i_1}y_{i_2}...y_{i_{r+1}} = \sum_{i_1<i_2<...<i_{r+1}}^{n} \frac{\sqrt[n]{I_{n,n}^{2(r+1)}}}{x_{i_1}x_{i_2}...x_{i_{r+1}}} = I_{n,n}^{\frac{2(r+1)}{n}-1} I_{n,n-r-1}.$$

Making use of (21), we will obtain

$$I_{n,n}^{\frac{2(r+1)}{n}-1} I_{n,n-r-1} = I_{n,n}^{1-\frac{2(n-r-1)}{n}} I_{n,n-r-1} = I_{n,r+1} \quad \to \quad J_{n,r+1} = I_{n,r+1}.$$

Continuing this procedure we come to the conclusion that the two sets $\{x_i\}$ and $\{y_i\}$ of real variables are related with $n$ equalities $J_{n,r} = I_{n,r}, r = 1, ..., n$. These equalities make coincident the roots $x_i$ and $y_i$ of both generating functions $F(x_i)$, $F(y_i)$

$$F(x_i) = \prod_{i=1}^{n}(x - x_i) = \sum_{r=0}^{n}(-1)^r I_{n,r} x^r = 0 = \sum_{r=0}^{n}(-1)^r J_{n,r} y^r = \prod_{i=1}^{n}(y - y_i) = F(y_i) ,$$

that leads to the coincidence of two sets $\{y_i\} = \{x_i\}$ up to the permutation $\widehat{\mathcal{D}}$ of the variables $y_i$ and $x_i$ in each set separately: $y_i = \widehat{\mathcal{D}} x_j$. This means in accordance with the definition (24) such $\lambda_{n*} = \sqrt[n]{I_{n,n}}$ always exists that for every $x_i$ one can point out such $x_j$, which satisfies $x_i \cdot x_j = \lambda_{n*}^2$ and $\lambda_{n*}$ is the common root of all $\mathsf{S}_{S_n}^{[m]}$- Eqns. ∎

At the end of this discussion we give simple examples for $\mathsf{S}_{S_3}^{[m]}$ and $\mathsf{S}_{S_4}^{[m]}$ polynomials, where the equivalence (21)⇔(22) is evident by factorization of invariant Eqns (21)

$$n = 3 : \quad I_{3,3} I_{3,1}^3 = I_{3,2}^3 \quad \to \quad (x_1^2 - x_2 x_3)(x_2^2 - x_3 x_1)(x_3^2 - x_1 x_2) = 0 .$$
$$n = 4 : \quad I_{4,4} I_{4,1}^2 = I_{4,3}^2 \quad \to \quad (x_1 x_2 - x_3 x_4)(x_1 x_3 - x_2 x_4)(x_1 x_4 - x_2 x_3) = 0 .$$

In general case a common solution $\lambda_n(x_i)$ of $\mathsf{S}_{S_n}^{\{m\}}$- Eqns for any $m$ does not exist. Now we will look for the special kind of such Eqns with non-negative coefficients $C_{s,l}$. This restriction makes it possible to find the universal bounds for $\lambda_n(x_i)$ that are valid for any $m$.

**Proposition 9**
*The unique positive root $\lambda_n(x_i)$ of $\mathsf{S}_{S_n}^{\{m\}}(\lambda_n, x_i)$- Eqn with $C_{s,l} \geq 0$ satisfies the following inequalities*

$$n \left( \sum_{i=1}^{n} \frac{1}{x_i} \right)^{-1} \leq \lambda_n(x_i) \leq \frac{1}{n} \sum_{i=1}^{n} x_i . \tag{25}$$

*Proof.* Making use of (14) it is sufficiently to prove only one of two inequalities (25), e.g. the right h.s. $\lambda_n \leq I_{n,1}/n$. Let us estimate the differences $\Gamma_{s,l} = T_\oplus^{s,l}\left(\genfrac{}{}{0pt}{}{\{m\}}{S_n}\right) - T_\ominus^{s,l}\left(\genfrac{}{}{0pt}{}{\{m\}}{S_n}\right)$

$$\Gamma_{s,l} = \prod_{r=1}^{n} I_{n,r}^{\alpha_r^{s,l}} \cdot \lambda_n^{mn-s} - I_{n,n}^{(m-\sum_{r=1}^{n} \alpha_r^{s,l})} \prod_{r=1}^{n} I_{n,n-r}^{\alpha_r^{s,l}} \cdot \lambda_n^s , \quad s = \sum_{r=1}^{n} r\, \alpha_r^{s,l} . \tag{26}$$

First we recall the inequalities [6] which relate the basic polynomials :

$$\left( \frac{I_{n,r}}{\binom{n}{r}} \right)^{\frac{1}{r}} \geq \left( \frac{I_{n,r+1}}{\binom{n}{r+1}} \right)^{\frac{1}{r+1}} , \quad \left( \frac{I_{n,n-r}}{\binom{n}{r} \cdot I_{n,n}} \right)^{\frac{1}{r}} \geq \left( \frac{I_{n,n-r-1}}{\binom{n}{r+1} \cdot I_{n,n}} \right)^{\frac{1}{r+1}} , \quad \text{i.e.} \tag{27}$$

$$\mathfrak{m}_a = \frac{1}{n} \cdot I_{n,1} \geq \left( 2 \frac{I_{n,2}}{n(n-1)} \right)^{\frac{1}{2}} \geq ... \geq \left( \frac{1}{n} \cdot I_{n,n-1} \right)^{\frac{1}{n-1}} \geq (I_{n,n})^{\frac{1}{n}} = \mathfrak{m}_g ,$$

$$\mathfrak{m}_g = (I_{n,n})^{\frac{1}{n}} \geq \left( \frac{n \cdot I_{n,n}}{I_{n,1}} \right)^{\frac{1}{n-1}} \geq ... \geq \left( \frac{n(n-1)}{2} \frac{I_{n,n}}{I_{n,n-2}} \right)^{\frac{1}{2}} \geq \frac{n \cdot I_{n,n}}{I_{n,n-1}} = \mathfrak{m}_h,$$

where $\mathfrak{m}_a, \mathfrak{m}_g, \mathfrak{m}_h$ denote the arithmetic, geometric and harmonic averages respectively.

Because of cumbersomeness of (26) we will hint the general result by evaluation the first successive inequalities making use of (27):

$$\Gamma_{0,1} = \lambda_n^{mn} - I_{n,n}^m \geq \lambda_n^{mn} - \left(n^{-1}I_{n,1}\right)^{mn}, \tag{28}$$

$$\Gamma_{1,1} = \lambda_n \cdot (I_{n,1}\lambda_n^{mn-2} - I_{n,n-1}I_{n,n}^{m-1}) \geq \lambda_n \cdot \left[I_{n,1}\lambda_n^{mn-2} - \frac{I_{n,1}^{n-1}}{n^{n-2}}\left(\frac{I_{n,1}}{n}\right)^{n(m-1)}\right] =$$

$$= \lambda_n I_{n,1} \cdot \left[\lambda_n^{mn-2} - \left(n^{-1}I_{n,1}\right)^{mn-2}\right],$$

$$\Gamma_{2,1} = \lambda_n^2 \cdot (I_{n,1}^2\lambda_n^{mn-4} - I_{n,n-1}^2 I_{n,n}^{m-2}) \geq \lambda_n^2 I_{n,1}^2 \cdot \left[\lambda_n^{mn-4} - \left(n^{-1}I_{n,1}\right)^{mn-4}\right],$$

$$\Gamma_{2,2} = \lambda_n^2 \cdot (I_{n,2}\lambda_n^{mn-4} - I_{n,n-2}I_{n,n}^{m-1}) \geq \lambda_n^2 I_{n,2} \cdot \left[\lambda_n^{mn-4} - \left(n^{-1}I_{n,1}\right)^{mn-4}\right],$$

$$\Gamma_{3,1} = \lambda_n^3 \cdot (I_{n,1}^3\lambda_n^{mn-6} - I_{n,n-1}^3 I_{n,n}^{m-3}) \geq \lambda_n^3 I_{n,1}^3 \cdot \left[\lambda_n^{mn-6} - \left(n^{-1}I_{n,1}\right)^{mn-6}\right].$$

Continuing this procedure we arrive at general result

$$\Gamma_{s,l} \geq \lambda_n^s \prod_{r=1}^{n} I_{n,r}^{\alpha_r^{s,l}} \cdot \left[\lambda_n^{mn-2s} - \left(n^{-1}I_{n,1}\right)^{mn-2s}\right]. \tag{29}$$

Making partial summation in (5) we obtain

$$C_{1,1}\Gamma_{1,1} \geq \lambda_n \, C_{1,1}I_{n,1} \left[\lambda_n^{mn-2} - \left(n^{-1}I_{n,1}\right)^{mn-2}\right],$$

$$\sum_{l=1}^{2} C_{2,l}\Gamma_{2,l} \geq \lambda_n^2 \left\{C_{2,1}I_{n,1}^2 + C_{2,2}I_{n,2}\right\} \left[\lambda_n^{mn-4} - \left(n^{-1}I_{n,1}\right)^{mn-4}\right],$$

$$\sum_{l=1}^{3} C_{3,l}\Gamma_{3,l} \geq \lambda_n^3 \left\{C_{3,1}I_{n,1}^3 + C_{3,2}I_{n,1}I_{n,2} + C_{3,3}I_{n,3}\right\} \left[\lambda_n^{mn-6} - \left(n^{-1}I_{n,1}\right)^{mn-6}\right], \; etc$$

where the curl brackets are always positive due to assumption $C_{s,l} \geq 0$. Using now the factorization $a^k - b^k = (a-b)\sum_{j=0}^{k-1} a^j \cdot b^{k-1-j}$ and substituting (29) into (5) we will get

$$\mathsf{S}_{S_n}^{\{m\}}(\lambda_n, x_i) = 0 \geq (\lambda_n - \frac{1}{n}I_{n,1}) \cdot \sum_{s=0}^{mn} \lambda_n^s \sum_{l=1}^{\mathcal{P}_n^m(s)} C_{s,l} \prod_{r=1}^{n} I_{n,r}^{\alpha_r^{s,l}} \sum_{j=0}^{mn-2s-1} \left(\frac{1}{n}I_{n,1}\right)^j \lambda_n^{mn-2s-1-j}.$$

The above inequality yields $\lambda_n(x_i) \leq \frac{1}{n} \cdot I_{n,1}$ that proves the Proposition 9. ∎

## 2.3 Monotonicity and universal bounds for $n = 3, 4$

In this Section we will show that the upper and lower bounds (25) for the non-negative algebraic solution $\lambda_n(x_i)$ of $\mathsf{S}_{S_n}^{\{m\}}$- Eqns with $C_{s,l} \geq 0$ can be enhanced if additional requirement of monotonicity of such bounds is assumed. We will consider two simple cases for $n = 3, 4$.

Let us rewrite the first two $\mathsf{S}_{S_3}^{\{m\}}$-equations in the following way

$$\mathsf{S}_{S_3}^{\{1\}} = \lambda_3^3 - I_{3,3} + A_3 I_{3,1}\lambda_3 \left(\lambda_3 - \frac{I_{3,2}}{I_{3,1}}\right), \quad A_3 \geq 0, \; D_3 \geq 0, \; B_3 \geq 0, \; C_3 \geq 0, \tag{30}$$

$$\mathsf{S}_{S_3}^{\{2\}} = \lambda_3^6 - I_{3,3}^2 + D_3 I_{3,1}\lambda_3 \left(\lambda_3^4 - \frac{I_{3,2}I_{3,3}}{I_{3,1}}\right) + B_3 I_{3,1}^2 \lambda_3^2 \left(\lambda_3^2 - \frac{I_{3,2}^2}{I_{3,1}^2}\right) + C_3 I_{3,2}\lambda_3^2 \left(\lambda_3^2 - \frac{I_{3,1}I_{3,3}}{I_{3,2}}\right).$$

Continuing this procedure we will get different expressions in brackets of generic type

$$\lambda_3^{a+3b} - \left(\frac{I_{3,2}}{I_{3,1}}\right)^a I_{3,3}^b, \quad a + 3b > 0,$$

where $a, b$ are integers. Indeed, according to (5) two counter-partners in $\mathsf{S}_{S_3}^{[m]}$-equation will be coupled as following

$$\lambda_3^{3m-k} I_{3,1}^u I_{3,2}^v I_{3,3}^w - \lambda_3^k I_{3,1}^v I_{3,2}^u I_{3,3}^{m-u-v-w} = \lambda_3^k I_{3,1}^u I_{3,2}^v I_{3,3}^w \left[\lambda_3^{3m-2k} - \left(\frac{I_{3,2}}{I_{3,1}}\right)^{u-v} I_{3,3}^{m-u-v-2w}\right]$$

$$= \lambda_3^k I_{3,1}^u I_{3,2}^v I_{3,3}^w \left[\lambda_3^{a+3b} - \left(\frac{I_{3,2}}{I_{3,1}}\right)^a I_{3,3}^b\right],$$

where $a = u - v$, $b = m - u - v - 2w$, $k = u + 2v + 3w$, $u, v, w \geq 0$. In order to find the upper $\Omega_3(x_i)$ and the lower $\omega_3(x_i)$ bounds we have to consider the following term

$$\Upsilon_3(p_3, x_i) = \left[I_{3,3}\left(\frac{I_{3,2}}{I_{3,1}}\right)^{p_3}\right]^{1/(3+p_3)}, \quad \text{where} \quad p_3 = \frac{a}{b}. \tag{31}$$

Let $\Omega_3(x_i) = max\{\Upsilon_3(p_3, x_i); -\infty \leq p_3 \leq \infty\}$, $\omega_3(x_i) = min\{\Upsilon_3(p_3, x_i); -\infty \leq p_3 \leq \infty\}$:

$$\Omega_3(x_i) = \Upsilon_3(p_3^+, x_i), \quad \omega_3(x_i) = \Upsilon_3(p_3^-, x_i).$$

Then we can rewrite Eqns (30) as following

$$\lambda^3 - \Omega_3^3(x_i) + A_3 I_{3,1} \lambda \left[\lambda - \Omega_3(x_i)\right] \leq 0 \leq \lambda^3 - \omega_3^3(x_i) + A_3 I_{3,1} \lambda \left[\lambda - \omega_3(x_i)\right],$$
$$\lambda^6 - \Omega_3^6(x_i) + D_3 I_{3,1} \lambda \left[\lambda^4 - \Omega_3^4(x_i)\right] + \{B_3 I_{3,1}^2 + C_3 I_{3,2}\}\lambda^2 \left[\lambda^2 - \Omega_3^2(x_i)\right] \leq 0,$$
$$\lambda^6 - \omega_3^6(x_i) + D_3 I_{3,1} \lambda \left[\lambda^4 - \omega_3^4(x_i)\right] + \{B_3 I_{3,1}^2 + C_3 I_{3,2}\}\lambda^2 \left[\lambda^2 - \omega_3^2(x_i)\right] \geq 0,$$

where the curl brackets are always positive according to (30). Using now the factorization $c^k - d^k$ by $c - d$ we arrive at the upper $\Omega_3(x_i)$ and the lower $\omega_3(x_i)$ bounds which govern the positive solutions of the first two $\mathsf{S}_{S_3}^{[m]}$-equations as well as for any other $m$. Now we will look for those $p_3$ which provide monotonicity of $\Upsilon_3(p_3, x_i)$. Let us find the first derivative

$$\frac{3+p_3}{\Upsilon_3} \cdot \frac{\partial \Upsilon_3}{\partial x_1} = \frac{1}{x_1} + p_3 \left(\frac{x_2 + x_3}{I_{3,2}} - \frac{1}{I_{3,1}}\right) = \frac{U_3}{x_1 I_{3,2} I_{3,1}}, \quad U_3 = I_{3,1} I_{3,2} + p_3 x_1 (x_2^2 + x_3^2 + x_2 x_3)$$

The monotonicity holds for $0 \leq p_3 \leq \infty$. One can show that this region might be enlarged up to $-1 \leq p_3 \leq \infty$. Indeed, assume $p_3 = -1 + \varepsilon_3$. Then

$$U_3 = (1 - \varepsilon_3)(x_1 I_{3,2} + x_2 x_3 I_{3,1}) + \varepsilon_3 I_{3,1} I_{3,2} \geq 0, \quad \text{if} \quad 0 \leq \varepsilon_3 \leq 1.$$

In contrast, if $\varepsilon_3 \leq 0$ we always can choose such values $x_3 \to 0$ and $x_1/x_2 + \varepsilon_3 \leq 0$ that

$$U_3 = (1 - \varepsilon_3) x_1^2 x_2 + \varepsilon_3 x_1 x_2 (x_1 + x_2) = x_1^2 x_2 \left[1 + \varepsilon_3 \frac{x_2}{x_1}\right] \leq 0.$$

Finally we get $p_3 \in [-1, \infty)$. In order to find the extremal $p_3^\pm$ which correspond to the upper and lower bounds respectively we should check the derivative over $p_3$

$$\frac{(3+p_3)^2}{\Upsilon_3} \cdot \frac{\partial \Upsilon_3(p_3, x_i)}{\partial p_3} = \ln\left[\left(\frac{I_{3,2}}{I_{3,1}}\right)^3 \frac{1}{I_{3,3}}\right] = \begin{cases} < 0, & x_i \in I_{3,3} I_{3,1}^3 \geq I_{3,2}^3, \\ = 0, & x_i \in I_{3,3} I_{3,1}^3 = I_{3,2}^3, \\ > 0, & x_i \in I_{3,3} I_{3,1}^3 \leq I_{3,2}^3 \end{cases}$$

The monotonicity of $\Upsilon_3$ in $p_3$ gives $p_3^{\pm} = -1, \infty$ and the bounds are the following

$$\omega_3(x_i) = min \left\{ \frac{I_{3,2}}{I_{3,1}}, \sqrt{I_{3,3} \cdot \frac{I_{3,1}}{I_{3,2}}} \right\}, \quad \Omega_3(x_i) = max \left\{ \frac{I_{3,2}}{I_{3,1}}, \sqrt{I_{3,3} \cdot \frac{I_{3,1}}{I_{3,2}}} \right\}. \quad (32)$$

The similar approach undertaken for $\mathsf{S}_{S_4}^{\{m\}}$- Eqn gives an admitting term

$$\Upsilon_4(p_4, x_i) = \left[ I_{4,4} \left( \frac{I_{4,3}}{I_{4,1}} \right)^{p_4} \right]^{1/(4+2p_4)}. \quad (33)$$

Let us find the first derivative

$$\frac{4+2p_4}{\Upsilon_4} \cdot \frac{\partial \Upsilon_4}{\partial x_1} = \frac{1}{x_1} + p_4 \left( \frac{x_2 x_3 + x_3 x_4 + x_4 x_2}{I_{4,3}} - \frac{1}{I_{4,1}} \right) = \frac{U_4}{x_1 I_{4,3} I_{4,1}},$$

$$U_4 = I_{4,1} I_{4,3} + p_4 x_1 \left[ x_2^2(x_3 + x_4) + x_3^2(x_4 + x_2) + x_4^2(x_2 + x_3) + 2x_2 x_3 x_4 \right],$$

that gives monotonicity for $0 \leq p_4 \leq \infty$. As in previous case $(n = 3)$ this region can be enlarged up to $-1 \leq p_4 \leq \infty$. Indeed, defining $p_4 = -1 + \varepsilon_4$, we have

$$U_4 = (1 - \varepsilon_4) \left[ x_1^2(x_2 x_3 + x_3 x_4 + x_4 x_2) + x_2 x_3 x_4(x_2 + x_3 + x_4) + 2x_2 x_3 x_4 \right] + \varepsilon_4 I_{4,1} I_{4,3} \geq 0,$$

if $0 \leq \varepsilon_4 \leq 1$. Taking $x_4 \to 0$ we obtain

$$U_4 = x_1^2 x_2 x_3 \left[ 1 + \varepsilon_4 \left( \frac{x_2}{x_1} + \frac{x_3}{x_1} \right) \right] \leq 0, \quad \text{if} \quad \frac{1}{\varepsilon_4} + \frac{x_2}{x_1} + \frac{x_3}{x_1} \leq 0.$$

Thus, we get finally $p_4 \in [-1, \infty)$. To find the extremal $p_4^{\pm}$ which correspond to the upper $\Omega_4(x_i)$ and lower $\omega_4(x_i)$ bounds respectively we should check the derivative over $p_4$

$$\frac{(4+2p_4)^2}{\Upsilon_4} \cdot \frac{\partial \Upsilon_4(p_4, x_i)}{\partial p_4} = \ln \left[ \left( \frac{I_{4,3}}{I_{4,1}} \right)^2 \frac{1}{I_{4,4}} \right] = \begin{cases} < 0, & x_i \in I_{4,4} I_{4,1}^2 \geq I_{4,3}^2, \\ = 0, & x_i \in I_{4,4} I_{4,1}^2 = I_{4,3}^2, \\ > 0, & x_i \in I_{4,4} I_{4,1}^2 \leq I_{4,3}^2 \end{cases}$$

The monotonicity of $\Upsilon_4$ in $p_4$ gives $p_4^{\pm} = -1, \infty$ and the bounds are the following

$$\omega_4(x_i) = min \left\{ \sqrt{\frac{I_{4,3}}{I_{4,1}}}, \sqrt{I_{4,4} \cdot \frac{I_{4,1}}{I_{4,3}}} \right\}, \quad \Omega_4(x_i) = max \left\{ \sqrt{\frac{I_{4,3}}{I_{4,1}}}, \sqrt{I_{4,4} \cdot \frac{I_{4,1}}{I_{4,3}}} \right\}. \quad (34)$$

## 3 Extension on the finite groups

In this Section we are going to discuss the existence of other finite groups $G$ distinguished from $S_n$ that makes it possible to build out the reciprocal $\mathsf{R}_G(\lambda_n, x_i)$ and skew-reciprocal $\mathsf{S}_G(\lambda_n, x_i)$ Eqns on the basis of homogeneous polynomial invariants $I_{d_k}(G)$ of degrees $d_k$. This question is motivated both from the physical and mathematical standpoints. The self-dual problem of effective isotropic conductivity in 2D 3-component regular and random checkerboards admits the tessellations of a plane with less rigid requirements of the permutation invariance upon the action of the cyclic group $\mathcal{Z}_3$ which also leads to skew-reciprocal Eqns [7]. Moreover, one can point out several finite groups where such property does exist and where it does not. We give two instructive examples concerned with finite 3D crystallographic groups.

*Cubic group $O_h$* [8]

$$\deg O_h = 3 \ , \ |O_h| = 48 \ , \ d_1 = 2 \ , \ d_2 = 4 \ , \ d_3 = 6 \ , \tag{35}$$
$$I_1(O_h) = x_1^2 + x_2^2 + x_3^2 \ , \ I_2(O_h) = x_1^2 x_2^2 + x_2^2 x_3^2 + x_3^2 x_1^2 \ , \ I_3(O_h) = x_1^2 x_2^2 x_3^2 \ ,$$
$$\mathsf{S}_{O_h}^{\{1\}}(\lambda, x_i) = a_0 \lambda^6 + a_1 I_1(O_h) \lambda^4 - a_1 I_2(O_h) \lambda^2 - a_0 I_3(O_h) = 0 \ , \ \mu\left\{\mathsf{R}_{O_h}^{\{1\}}\right\} = \mu\left\{\mathsf{S}_{O_h}^{\{1\}}\right\} = 2 \ .$$

*Axial group $C_n$* [8]

$$\deg C_n = 2 \ , \ |C_n| = n \ , \ d_1 = 2 \ , \ d_2 = n \ ,$$
$$I_1(C_n) = x_1^2 + x_2^2 \ , \ I_2(C_n) = x^n \cos n\phi \ , \ I_3(C_n) = x^n \sin n\phi \ ,$$
$$x_1 = x \cos \phi \ , \ x_2 = x \sin \phi \ , \ I_2^2(C_n) + I_3^2(C_n) = I_1^n(C_n) \ , \tag{36}$$

e.g.
$$I_1(C_3) = x_1^2 + x_2^2 \ , \ I_2(C_3) = x_1^3 - 3 \, x_1 \, x_2^2 \ , \ I_3(C_3) = x_2^3 - 3 \, x_2 \, x_1^2 \ .$$

The basic invariants of $C_n$ are related with syzygy of the 1-st kind (36) and do not give rise to the skew-reciprocal Eqn when $n \neq 2, 4$ that will be proven few lines later. This indicates a strong relationship between the existence of the colour tessellation with certain finite group $G$ that belongs to the family $\mathfrak{S}_W$ of finite groups, which can be matched with conformal group $\mathsf{W}$ and produce $\mathsf{R}_G(\lambda, x_i)$ and $\mathsf{S}_G(\lambda, x_i)$ Eqns. In the rest of the Section we will find the necessary condition of such matching.

**Lemma 1**
*If $\lambda(x_i)$ is a root of $\mathsf{S}_G(\lambda, x_i)$- Eqn then necessarily*

$$d_k + d_{n-k} = d_n \ , \ I_{d_n}(x_i) \cdot I_{d_k}\left(\frac{1}{x_i}\right) = I_{d_{n-k}}(x_i) \quad \rightarrow \quad I_{d_n}(x_i) \cdot I_{d_n}\left(\frac{1}{x_i}\right) = 1 \tag{37}$$

*and $I_{d_n}(x_i)$ is a monomial.*

*Proof.* Suppose a finite group $G$ acts on the $n$ dimensional vector space $V_n$ over the complex numbers $\mathbf{C}$, and a tuple $\{x_1, x_2, ..., x_n\}$ denotes a basis for $V_n$. There exist exactly $n$ algebraically independent invariants $I_{d_k}(x_i)$. A basis of algebra of invariants of $G$ is generated by not more than $\binom{|G|+n}{n}$ homogeneous invariants of degrees $d_k$ not exceeding the order $|G|$ of group. Let us write both $\mathsf{S}_G(\lambda, x_i)$ and $\mathsf{S}_G(\lambda^{-1}, x_i^{-1})$ monic polynomials

$$\mathsf{S}_G(\lambda, x_i) = \lambda^{d_n} + a_1 I_{d_1} \lambda^{d_n - d_1} + a_2 I_{d_2} \lambda^{d_n - d_2} + \ ... \ - a_2 I_{d_{n-2}} \lambda^{d_n - d_{n-2}} - a_1 I_{d_{n-1}} \lambda^{d_n - d_{n-1}} - I_{d_n} \ ,$$

$$\mathsf{S}_G\left(\frac{1}{\lambda}, \frac{1}{x_i}\right) = \lambda^{d_n} + a_1 \frac{J_{d_n}}{J_{d_{n-1}}} \lambda^{d_n - 1} + a_2 \frac{J_{d_n}}{J_{d_{n-2}}} \lambda^{d_n - 2} + \ ... \ - a_2 \frac{J_{d_n}}{J_{d_2}} \lambda^{d_2} - a_1 \frac{J_{d_n}}{J_{d_1}} \lambda^{d_1} - J_{d_n} \ ,$$

where $J_{d_k}(x_i) = I_{d_k}^{-1}\left(\frac{1}{x_i}\right)$ is rational function in $x_i$. If $G \in \mathfrak{S}_W$ then both polynomials $\mathsf{S}_G(\lambda, x_i)$ and $\mathsf{S}_G\left(\frac{1}{\lambda}, \frac{1}{x_i}\right)$ coincide and we arrive at (37). The monomiality of $I_{d_n}(x_i)$ is trivial and can be proven by opposite assumption. ∎

The degrees $d_k$ are uniquely determined, while the polynomials $I_{d_k}(x_i)$ themselves are not unique: they depend on the basis $\{x_1, x_2, ..., x_n\}$ chosen for $V_n$ and may be replaced by their algebraic combinations. This makes our Lemma 1 less effective. Nevertheless we are able to advance in its usage with help of the following auxiliary

**Lemma 2**
*Let $F(x_1, x_2)$ be a monic homogeneous form of degree n. If there exists a non-degenerated linear operator $\widehat{A}$, $(x_1, x_2) \xrightarrow{\widehat{A}} (y_1, y_2)$, that transforms $F(x_1, x_2)$ into monomial $y_1^{\alpha_1} y_2^{\alpha_2}$, $\alpha_1 + \alpha_2 = n$, then F has exactly two degenerated roots of degrees $\alpha_1$ and $\alpha_2$.*

<u>Proof.</u>   Let us decompose $F$ uniquely into a product of linear forms over the complex numbers **C** and apply a linear transformation $\widehat{A}$

$$F(x_1, x_2) = \prod_{i=1}^{n}(x_1 - c_i x_2)^{\nu_i} \xrightarrow{\widehat{A}} \prod_{i=1}^{n}(b'_i y_1 - c'_i y_2)^{\nu_i} , \qquad (38)$$

where $\nu_i$ are the degrees of degeneration. If the right h.s. of (38) is monomial in variables $y_1, y_2$ then there must be only two multipliers which contain zero coefficients, viz, $b'_2 = c'_1 = 0$ that yields

$$\nu_1 = \alpha_1 , \quad \nu_2 = \alpha_2 , \quad \nu_k = 0 , \ k \geq 3 .$$

By the other hand the existence of only two degenerated roots of $F(x_1, x_2)$ is sufficient to represent it as a monomial. ∎

Hence, the axial group $C_n, n \neq 2, 4$ does not give rise to the skew-reciprocal Eqn since both invariants

$$I_2(C_n) = \sum_{j=0}^{2j \leq n}(-1)^j \binom{n}{2j} x_1^{n-2j} x_2^{2j} , \quad I_3(C_n) = \sum_{j=0}^{2j \leq n-1}(-1)^j \binom{n}{2j+1} x_2^{n-2j-1} x_1^{2j+1}$$

have more than two real roots that follows from the Descartes' rule of signs interchange in the sequence of coefficients for the real algebraic equations. The cases $n = 2, 4$ are both exceptional: $I_3(C_2) = 2\,x_1\,x_2$, $I_1^2(C_4) - I_2(C_4) = 8\,x_1^2 x_2^2$.

We will consider now some classical groups $G \subset S_n$ where the existence of skew-reciprocal Eqns can be proved or disproved using Lemma 1.

## 3.1  Cyclic $\mathcal{Z}_n$ and alternating $\mathcal{A}_n$ groups

The cyclic group $\mathcal{Z}_n$ is an Abelian group and its invariants are generated by elementary symmetric functions $I_{n,r}$ (4) together with a product

$$I_{n+1}(\mathcal{Z}_n) = \prod_{i}^{n}(x_i - x_{i+1}) , \quad x_{n+1} = x_1 ,$$

whose square is a polynomial in the invariants $I_{n,r}$.

The alternating group $\mathcal{A}_n$ is a subgroup of index 2 of the full permutation group $S_n$. Its invariants are generated by elementary symmetric functions $I_{n,r}$ (4) together with the van der Monde discriminant [2]

$$I_{n+1}(\mathcal{A}_n) = \det \begin{pmatrix} 1 & x_1 & x_1^2 & \dots & x_1^{n-1} \\ 1 & x_2 & x_2^2 & \dots & x_2^{n-1} \\ . & . & . & \dots & . \\ 1 & x_n & x_n^2 & \dots & x_n^{n-1} \end{pmatrix} = \prod_{i<j}^{n}(x_i - x_j) ,$$

whose square is also a polynomial in the invariants $I_{n,r}$.

Their additional invariants $I_{n+1}(\mathcal{A}_n)$ and $I_{n+1}(\mathcal{Z}_n)$ are algebraically dependent and cannot be incorporated into $\mathsf{S}_{S_n}^{\{m\}}$- equation for $\lambda$ because that would violate the duality requirement. Thus, we are lead back to the $\mathsf{S}_{S_n}^{\{m\}}$-equation (3), which is dictated not only by the strong requirement of the full permutation invariance $S_n$, but even by the milder requirements of alternating $\mathcal{A}_n$ or cyclic $\mathcal{Z}_n$ permutation invariance.

## 3.2 Coxeter groups $\mathsf{B_n}, \mathsf{D_n}, \mathsf{I_{2,m}}, \mathsf{G_2}, \mathsf{H_3}$

The detailed description of invariant polynomials for Coxeter groups and their discussion is given in [9]. The maximal degree $d_n(G)$ is usually referred to as a Coxeter number.

$\mathsf{B_n}$ , $d_r(\mathsf{B_n}) = 2, 4, 6, ..., 2n$ .

$$I_r(\mathsf{B_n}) = \sum_{i_1 < i_2 < ... < i_r}^{n} x_{i_1}^2 x_{i_2}^2 ... x_{i_r}^2 . \tag{39}$$

$\mathsf{D_n}$ , $d_r(\mathsf{D_n}) = 2, 4, 6, ..., 2(n-1), n$ , $n \geq 4$ .

$$I_r(\mathsf{D_n}) = \sum_{i_1 < i_2 < ... < i_r}^{n} x_{i_1}^2 x_{i_2}^2 ... x_{i_r}^2 , \ r \leq n-1 , \quad I_n(\mathsf{D_n}) = \prod_{i=1}^{n} x_i . \tag{40}$$

The reflection groups $\mathsf{B_n}$ and $\mathsf{D_n}$ give rise to the skew-reciprocal Eqns since their action on a set of $n$ real variables $\{x_1, x_2, ..., x_n\}$ is similar to the action of symmetric group $S_n$ on the set of $n$ non-negative variables $\{x_1^2, x_2^2, ..., x_n^2\}$.

$\mathsf{I_{2,m}}$, $m \neq 2, 4$ , $d_r(\mathsf{I_{2,m}}) = 2, m$ , $\mathsf{G_2} = \mathsf{I_{2,6}}$.

$$I_1(\mathsf{I_{2,m}}) = x_1^2 + x_2^2 , \quad I_2(\mathsf{I_{2,m}}) = \sum_{j=1}^{m} \left( x_1 \cos \frac{2\pi j}{m} + x_2 \sin \frac{2\pi j}{m} \right)^m . \tag{41}$$

The dihedral group $\mathsf{I_{2,m}}$ produces the invariant $I_2(\mathsf{I_{2,m}})$ which has more than two roots in the complex plane; due to Lemma 2 it does not satisfy Lemma 1 and therefore does not lead to the skew-reciprocal Eqns. The exceptional cases are $\mathsf{I_{2,2}} = \mathsf{D_2}$, $\mathsf{I_{2,4}} = \mathsf{B_2}$.

$\mathsf{H_3} = Y_h$ , $d_r(\mathsf{H_3}) = 2, 6, 10$ .

Three basic invariants for icosahedral group $Y_h$ were derived in [10]

$$I_1(\mathsf{H_3}) = x_1^2 + x_2^2 + x_3^2 , \ I_2(\mathsf{H_3}) = (\tau^2 \, x_1^2 - x_2^2)(\tau^2 \, x_2^2 - x_3^2)(\tau^2 \, x_3^2 - x_1^2) , \ \tau = \frac{1 + \sqrt{5}}{2} ,$$
$$I_3(\mathsf{H_3}) = (x_1^4 + x_2^4 + x_3^4)(\tau^4 \, x_1^2 - x_2^2)(\tau^4 \, x_2^2 - x_3^2)(\tau^4 \, x_3^2 - x_1^2) . \tag{42}$$

A decomposable representation (42) of $I_2(\mathsf{H_3})$ and $I_3(\mathsf{H_3})$ makes it possible to prove an absence of linear transformation which reduces them to the monomials $y_1^{\alpha_1} y_2^{\alpha_2} y_3^{\alpha_3}$.

## 3.3 Symmetric groups' product $S_{n_1} \times S_{n_2}$

The permutation invariance with a direct product of groups $G_1 \times G_2$ appears in the colour tessellation of 2D $n$-component checkerboard [4] when one can decompose it in two lattices $L_1, L_2$ and every lattice is coloured randomly with $n_1$ and $n_2$ colours respectively, $n =$ $n_1 + n_2$. In the case of equipartition of the colours in every lattice separately we come to $G = S_{n_1} \times S_{n_2} \subset S_{n_1+n_2}$. Thus, a definition of basic polynomial invariants is the following

$$V = V_{n_1} \oplus V_{n_2} , \quad \deg(S_{n_1} \times S_{n_2}) = n_1 + n_2 , \quad |S_{n_1} \times S_{n_2}| = n_1! \cdot n_2!$$

$$X_{n_1,r} = \sum_{i_1 < i_2 < ... < i_r}^{n_1} x_{i_1} x_{i_2} ... x_{i_r} , \quad Y_{n_2,r} = \sum_{i_1 < i_2 < ... < i_r}^{n_2} y_{i_1} y_{i_2} ... y_{i_r}$$

and corresponding $\mathsf{S}^{\{m_1,m_2\}}_{S_{n_1} \times S_{n_2}}(\lambda, x_i, y_j)$- Eqn reads

$$\begin{aligned}
\mathsf{S}^{\{m_1,m_2\}}_{S_{n_1} \times S_{n_2}} &= \lambda^{m_1 n_1 + m_2 n_2} + [B_{1,1} X_{n_1,1} + C_{1,1} Y_{n_2,1}] \lambda^{m_1 n_1 + m_2 n_2 - 1} + \big[B_{2,1} X^2_{n_1,1} + \\
&\quad B_{2,2} X_{n_1,2} + D_{2,1} X_{n_1,1} Y_{n_2,1} + C_{2,1} Y^2_{n_2,1} + C_{2,2} Y_{n_2,2}\big] \lambda^{m_1 n_1 + m_2 n_2 - 2} + \\
&\quad \dotsb - \\
&\quad \big[(B_{2,1} X^2_{n_1,n_1-1} + B_{2,2} X_{n_1,n_1-2} X_{n_1,n_1}) Y^{m_2}_{n_2,n_2} + \\
&\quad D_{2,1} X_{n_1,n_1-1} Y_{n_2,n_2-1} X^{m_1-1}_{n_1,n_1} Y^{m_2-1}_{n_2,n_2} + \\
&\quad (C_{2,1} Y^2_{n_2,n_2-1} + C_{2,2} Y_{n_2,n_2-2} Y_{n_2,n_2}) X^{m_1}_{n_1,n_1}\big] \lambda^2 - \big[(B_{1,1} X_{n_1,n_1-1} X^{m_1-1}_{n_1,n_1} Y^{m_2}_{n_2,n_2} + \\
&\quad C_{1,1} Y_{n_2,n_2-1} Y^{m_2-1}_{n_2,n_2} X^{m_1}_{n_1,n_1}\big] \lambda - X^{m_1}_{n_1,n_1} Y^{m_2}_{n_2,n_2} = 0 .
\end{aligned} \quad (43)$$

Here we come to generalization of conformal partition function $\mathcal{P}^m_n(s)$ : every coefficient at the term $\lambda^{m_1 n_1 + m_2 n_2 - s}$ in (43) is build out of independent polynomial invariants $X_{n_1,r}$, $Y_{n_2,r}$ and their number $\mathcal{P}\begin{bmatrix}m_1,m_2\\n_1,n_2\end{bmatrix}(s)$ is related to usual conformal partitions (see Section 4.1)

$$\mathcal{P}\begin{bmatrix}m_1,m_2\\n_1,n_2\end{bmatrix}(s) = \sum_{s_1=0}^{s} \mathcal{P}^{m_1}_{n_1}(s_1) \cdot \mathcal{P}^{m_2}_{n_2}(s - s_1) , \quad \mathcal{P}\begin{bmatrix}m_1,m_2\\n_1,n_2\end{bmatrix}(s) = \mathcal{P}\begin{bmatrix}m_2,m_1\\n_2,n_1\end{bmatrix}(s) . \quad (44)$$

There is one more simplification of (43) due to the additional freedom of permutation of lattices $L_1$ and $L_2$ that diminishes essentially the unimodality indices in (43). Following Jahn notation [11] we mark the action of such group as $[S_{n_1} \times S_{n_2}]$ and write

$$\begin{aligned}
\mathsf{S}^{\{m_1,m_2\}}_{[S_{n_1} \times S_{n_2}]} &= \lambda^{m_1 n_1 + m_2 n_2} + B_{1,1}[X_{n_1,1} + Y_{n_2,1}] \lambda^{m_1 n_1 + m_2 n_2 - 1} + \\
&\quad \big[B_{2,1}(X^2_{n_1,1} + Y^2_{n_2,1}) + B_{2,2}(X_{n_1,2} + Y_{n_2,2}) + D_{2,1} X_{n_1,1} Y_{n_2,1}\big] \lambda^{m_1 n_1 + m_2 n_2 - 2} + \\
&\quad \dotsb - \\
&\quad [B_{2,1}(X^2_{n_1,n_1-1} Y^{m_2}_{n_2,n_2} + Y^2_{n_2,n_2-1} X^{m_1}_{n_1,n_1}) + D_{2,1} X_{n_1,n_1-1} Y_{n_2,n_2-1} X^{m_1-1}_{n_1,n_1} Y^{m_2-1}_{n_2,n_2} + \\
&\quad B_{2,2}(X_{n_1,n_1-2} X_{n_1,n_1} Y^{m_2}_{n_2,n_2} + Y_{n_2,n_2-2} Y_{n_2,n_2} X^{m_1}_{n_1,n_1})] \lambda^2 - \\
&\quad B_{1,1}[X_{n_1,n_1-1} X^{m_1-1}_{n_1,n_1} Y^{m_2}_{n_2,n_2} + Y_{n_2,n_2-1} Y^{m_2-1}_{n_2,n_2} X^{m_1}_{n_1,n_1}] \lambda - X^{m_1}_{n_1,n_1} Y^{m_2}_{n_2,n_2} = 0 .
\end{aligned} \quad (45)$$

Simple comparison of (3), (43) and (45) gives

$$\mu\left\{\mathsf{S}^{\{m\}}_{S_{n_1+n_2}}\right\} < \mu\left\{\mathsf{S}^{\{m_1,m_2\}}_{[S_{n_1} \times S_{n_2}]}\right\} < \mu\left\{\mathsf{S}^{\{m_1,m_2\}}_{S_{n_1} \times S_{n_2}}\right\} , \quad m = \left[\frac{m_1 n_1 + m_2 n_2}{n_1 + n_2}\right] .$$

Similar inequalities are also true for reciprocal unimodality indices.

The further extension on a multiple direct product of the symmetric groups $S_{n_i}$ and symmetrized action of this product is obvious.

# 4 Gaussian generating function and conformal partitions

Let us consider now $\mathcal{P}_n^m(s)$ as an arithmetic function. The corresponding Diophantine equation $\sum_{r=1}^n r\, x_r^j = s$ follows from the similar one for $\mathcal{W}_n(s)$. However the admitted invariants, which contribute to term $\lambda_n^{mn-s}$ must satisfy additional restriction $\sum_{r=1}^n x_r^j \leq m$. That yields immediately two important facts

$$\mathcal{P}_n^m(s) \leq \mathcal{W}_n(s) \ , \quad \lim_{m \to \infty} \mathcal{P}_n^m(s) = \mathcal{W}_n(s) \ . \tag{46}$$

We will give more rigorous proof for (46) later in this Section.

The contraction of the conformal partition problem to the Diophantine system established above provides a big advantage. This system is strongly related to the Gaussian generating function and corresponding Cayley-Sylvester theorem [12].

**Theorem 1** (Cayley-Sylvester)
*Let $\mathcal{P}_n^m(s)$ is a number of the positive integer solutions $X = \{x_r\}$ of the system of one Diophantine equation and one Diophantine inequality*

$$\sum_{r=1}^n r\, x_r = s \ , \quad \sum_{r=1}^n x_r \leq m \ . \tag{47}$$

*1. Then $\mathcal{P}_n^m(s)$, which is a number of partitions of $s$ into at most $m$ parts, each $\leq n$, is generated by the Gaussian polynomial $G(n, m; t)$ of the finite order $nm$*

$$G(n, m; t) = \frac{\prod_{i=1}^{n+m}(1-t^i)}{\prod_{j=1}^{n}(1-t^j) \cdot \prod_{k=1}^{m}(1-t^k)} = \sum_{s \geq 0}^{nm} \mathcal{P}_n^m(s) \cdot t^s \ . \tag{48}$$

*2. $\mathcal{P}_n^m(s)$ has the following properties:*

$$\begin{aligned}
&\mathcal{P}_n^m(s) = 0 \quad \text{if} \ s > nm \ , \quad \mathcal{P}_n^m(0) = \mathcal{P}_n^m(nm) = 1 \ , \\
&\mathcal{P}_n^m(s) = \mathcal{P}_m^n(s) = \mathcal{P}_n^m(nm-s) \ , \quad \mathcal{P}_n^m\left(\frac{mn}{2}-s\right) = \mathcal{P}_n^m\left(\frac{mn}{2}+s\right) \ , \\
&\mathcal{P}_n^m(s) - \mathcal{P}_n^m(s-1) \geq 0 \quad \text{for} \quad 0 < s \leq \frac{nm}{2} \ .
\end{aligned} \tag{49}$$

Standard way to derive the expressions for $\mathcal{P}_n^m(s)$ is to study the recursive relations generated by $G(n, m; t)$. Such procedure was chosen in [3] for evaluation of Sylvester waves $\mathcal{W}_n(s)$.

Here we are going to clarify the basic properties of $\mathcal{P}_n^m(s)$ exploiting another approach, namely, the relationship between Molien (1) and Gaussian (48) generating functions

$$G(n, m; t) = \frac{M(n, t) \cdot M(m, t)}{M(n+m, t)} \ . \tag{50}$$

This gives rise to the following linear system

$$\Delta = \sum_{s=0}^{g} [\,\mathcal{P}_n^m(s) \cdot \mathcal{W}_{n+m}(g-s) - \mathcal{W}_n(s) \cdot \mathcal{W}_m(g-s)\,] = 0 \ , \quad \mathcal{W}_n(0) = \mathcal{W}_n(1) = 1 \ . \tag{51}$$

The Eqn (51) represents the linear convolution equations with a triangular Toeplitz matrix. It can be solved using the formulas [13] that express the inverse of the Toeplitz matrix. We will give another representation of a general solution of (51) which can be found due to triangularity of the matrix (see Appendix A)

$$\mathcal{P}_n^m(g) = \sum_{r=0}^{g-1} \left\{ \sum_{s=0}^{g-r} \mathcal{W}_n(s) \cdot \mathcal{W}_m(g-r-s) - \mathcal{W}_{n+m}(g-r) \right\} \cdot \Phi_r(-\mathcal{W}_{n+m}) \,, \quad \mathcal{P}_n^m(0) = 1 \,, \tag{52}$$

where polynomials $\Phi_r(-\mathcal{W}_{n+m})$ are the following

$$\Phi_r(-\mathcal{W}_{n+m}) = \sum_{q=1}^{r} q! \sum_{\{q_l\}}^{q=const} \frac{1}{q_1! \cdot q_2! \cdot \ldots \cdot q_r!} \prod_{l=1}^{r} (-1)^{q_l} \mathcal{W}_{n+m}^{q_l}(l) \,, \quad \sum_{l=1}^{r} q_l \cdot l = r \,, \quad \sum_{l=1}^{r} q_l = q \,.$$

The formula (52) solves basically the problem in closed form. Despite of its cumbersomeness this formula makes it possible to establish one important fact about the periodic properties of conformal partition function $\mathcal{P}_n^m(s)$. Indeed, according to [3] the first multiplier in (52) contains the periodic functions of all successive periods

$$\tau\{\mathcal{W}_n(s) \cdot \mathcal{W}_m(g-r-s)\} = 1, 2, ..., m \cdot n \,, \quad \tau\{\mathcal{W}_{n+m}(g-r)\} = 1, 2, ..., m+n \,,$$

while the second multiplier contains the periodic functions of periods $\tau\{\Phi_r(-\mathcal{W}_{n+m})\} = 1, 2, ..., m+n$. Hence due to the trivial inequality $m \cdot n \geq m + n$, $m, n \geq 2$ the conformal partition function $\mathcal{P}_n^m(s)$ contains the periodic functions at most of periods

$$\tau\{\mathcal{P}_n^m(s)\} = \begin{cases} 1, 2, ..., m \cdot n \,, & \text{if } m, n \geq 2 \,, \\ 1, 2, ..., m+1 \,, & \text{if } n = 1 \,, \\ 1, 2, ..., n+1 \,, & \text{if } m = 1 \,. \end{cases} \tag{53}$$

It turns out that the solution (52) for small $s$ can be essentially simplified. Assuming hereafter $n \leq m$ we will prove

**Proposition 10**

$$\mathcal{P}_n^m(s) = \begin{cases} \mathsf{P}(s) \,, & 0 \leq s \leq n \,, \\ \mathcal{W}_n(s) \,, & n \leq s \leq m \,, \\ \mathcal{V}_n^m(s) \,, & m \leq s \leq m+n \,, \end{cases} \tag{54}$$

where $\mathsf{P}(s)$ is unrestricted partition function and $\mathcal{V}_n^m(s) = \mathcal{W}_n(s) + \mathcal{W}_m(s) - \mathcal{W}_{m+n}(s)$.

*Proof.* Indeed, the first two lines in (54) follow from (51) by use of simple identities for restricted and unrestricted partition functions :

$$\begin{aligned} 0 \leq s \leq n &\quad \to \quad \mathsf{P}(s) = \mathcal{W}_{m+n}(s) = \mathcal{W}_m(s) = \mathcal{W}_n(s) \,, \\ n \leq s \leq m &\quad \to \quad \mathsf{P}(s) = \mathcal{W}_{m+n}(s) = \mathcal{W}_m(s) \,. \end{aligned} \tag{55}$$

Let us return now to (51) for $g = m + r$ and decompose it in the following way

$$\begin{aligned} \Delta &= \Delta_1 + \Delta_2 = 0 \,, \quad \Delta_1 = \sum_{s=0}^{r} F_1(s) \,, \quad \Delta_2 = \begin{cases} \sum_{s=r+1}^{k} F_1(s), & \text{if } m+r = 2k+1 \,, \\ \sum_{s=r+1}^{k-1} F_1(s) + F_2(k), & \text{if } m+r = 2k \,, \end{cases} \\ F_1(s) &= \mathcal{P}_n^m(s) \cdot \mathcal{W}_{n+m}(m+r-s) - \mathcal{W}_n(s) \cdot \mathcal{W}_m(m+r-s) + \\ &\quad \mathcal{P}_n^m(m+r-s) \cdot \mathcal{W}_{n+m}(s) - \mathcal{W}_n(m+r-s) \cdot \mathcal{W}_m(s) \,, \\ F_2(s) &= \mathcal{P}_n^m(s) \cdot \mathcal{W}_{n+m}(s) - \mathcal{W}_n(s) \cdot \mathcal{W}_m(s) \,. \end{aligned}$$

Suppose $r \leq n$ that leads to $k \leq m$ and $F_2(k) = 0$. The term $\Delta_2$ vanishes due to (54, 55):

$$s \leq m \to \begin{cases} \mathcal{P}_n^m(s) = \mathcal{W}_n(s) \\ \mathcal{W}_{m+n}(s) = \mathcal{W}_m(s) \end{cases}, \quad r \leq s \to \begin{cases} \mathcal{W}_{n+m}(m+r-s) = \mathcal{W}_m(m+r-s) \\ \mathcal{P}_n^m(m+r-s) = \mathcal{W}_n(m+r-s) \end{cases},$$

while the term $\Delta_1$ for $s \leq n$ converts into

$$\Delta_1 = \sum_{s=0}^{r} \mathsf{P}(s) \cdot [\mathcal{W}_{n+m}(m+r-s) - \mathcal{W}_m(m+r-s) + \mathcal{P}_n^m(m+r-s) - \mathcal{W}_n(m+r-s)],$$

that leads for $s \leq m + n$ to the third line in (54). ∎

Now we will give a brief proof of the limit in (46).

**Corollary 1**

$$\lim_{m \to \infty} \mathcal{P}_n^m(s) = \mathcal{W}_n(s), \quad \lim_{n \to \infty} \mathcal{P}_n^m(s) = \mathcal{W}_m(s). \tag{56}$$

*Proof.* The first limit follows from the 2-nd line in (54). The second limit follows from the exchange symmetry of indices (49). ∎

The further extension of (54) on $s > m+n$ loses its generosity and is more sophisticated. Here we are going to prove the following

**Proposition 11**

$$\mathcal{P}_n^m(m+n+k) = \mathcal{V}_n^m(m+n+k) + \mathcal{D}_n^m(k), \quad 0 \leq k \leq m-n, \quad k \leq n+1, \tag{57}$$

where $\mathcal{D}_n^m(k) = \mathcal{A}_n^m(k) - \mathcal{B}_n^m(k)$, $\mathcal{D}_n^m(0) = 0$ and

$$\mathcal{A}_n^m(k) = \sum_{k_1,k_2=0}^{k_1+k_2=k} [\mathsf{P}(m+k_1) - \mathcal{W}_m(m+k_1)] \cdot [\mathsf{P}(n+k_2) - \mathcal{W}_n(n+k_2)],$$

$$\mathcal{B}_n^m(k) = \sum_{r=0}^{k-1} \mathsf{P}(k-r) \cdot \mathcal{D}_n^m(r), \quad \mathcal{D}_n^m(r) = \mathcal{D}_m^n(r), \quad i.e.: \tag{58}$$

$$\mathcal{D}_n^m(1) = 0, \quad \mathcal{D}_n^m(2) = \mathcal{A}_n^m(2), \quad \mathcal{D}_n^m(3) = \mathcal{A}_n^m(3) - \mathsf{P}(1)\mathcal{D}_n^m(2), \text{ etc.}$$

*Proof.* We will prove (57) by induction. First, it is easy to check from (51) that (57) gives correct expression for $k = 1$. We assume that it is also true for $\mathcal{P}_n^m(m+n+k-r)$ for all $r \geq 1$. Let us decompose Eqn (51) in several parts

$$\mathcal{P}_n^m(m+n+k) = \mathcal{V}_n^m(m+n+k) + \mathsf{K}_1 + \mathsf{K}_2 + \mathsf{K}_3 + \mathsf{K}_4, \tag{59}$$

$$\mathsf{K}_1 = \sum_{r=1}^{k-2} F_3(r), \quad \mathsf{K}_2 = \sum_{r=k-1}^{n} F_3(r), \quad \mathsf{K}_3 = \sum_{r=n+1}^{n+k-1} F_3(r), \quad \mathsf{K}_4 = \sum_{r=n+k}^{[\frac{m+n+k}{2}]} F_3(r),$$

$$F_3(r) = \mathcal{W}_m(r) \cdot \mathcal{W}_n(m+n+k-r) - \mathcal{P}_n^m(m+n+k-r) \cdot \mathcal{W}_{m+n}(r) +$$
$$\mathcal{W}_m(m+n+k-r) \cdot \mathcal{W}_n(r) - \mathcal{P}_n^m(r) \cdot \mathcal{W}_{m+n}(m+n+k-r)$$

and find all contributions $\mathsf{K}_i$ in four different regions.

$\underline{1 \leq r \leq k-2 < n, \ m+n+2 \leq m+n+k-r \leq m+n+k-1}$

$$\begin{cases} \mathcal{W}_m(r) = \mathcal{W}_n(r) = \mathcal{W}_{m+n}(r) = \mathcal{P}_n^m(r) = \mathsf{P}(r) \\ \mathcal{P}_n^m(m+n+k-r) - \mathcal{V}_n^m(m+n+k-r) = \mathcal{D}_n^m(k-r) \end{cases} \longrightarrow$$

$$\mathsf{K}_1 = -\sum_{r=1}^{k-2} \mathsf{P}(r) \cdot \mathcal{D}_n^m(k-r) = -\sum_{r=0}^{k-1} \mathsf{P}(k-r) \cdot \mathcal{D}_n^m(r) \,, \ \text{since} \ \mathcal{D}_n^m(1) = 0 \,.$$

$\underline{k-1 \leq r \leq n, \ m+k \leq m+n+k-r \leq m+n+1}$

$$\begin{cases} \mathcal{W}_m(r) = \mathcal{W}_n(r) = \mathcal{W}_{m+n}(r) = \mathcal{P}_n^m(r) \\ \mathcal{P}_n^m(m+n+k-r) = \mathcal{V}_n^m(m+n+k-r) \end{cases} \longrightarrow \mathsf{K}_2 = 0 \,.$$

$\underline{n+1 \leq r \leq n+k-1 < m, \ m+1 \leq m+n+k-r \leq m+k-1 < m+n}$

$$\begin{cases} \mathcal{P}_n^m(r) = \mathcal{W}_n(r) \\ \mathcal{W}_{m+n}(r) = \mathcal{W}_m(r) \end{cases} \text{and} \begin{cases} \mathcal{W}_n(s) - \mathcal{P}_n^m(s) = \mathcal{W}_{m+n}(s) - \mathcal{W}_m(s) \\ m+1 \leq s = m+n+k-r \leq m+n \end{cases} \longrightarrow \quad (60)$$

$$\mathsf{K}_3 = \sum_{r=n+1}^{n+k-1} [\mathcal{W}_{m+n}(m+n+k-r) - \mathcal{W}_m(m+n+k-r)] \cdot [\mathcal{W}_{m+n}(r) - \mathcal{W}_n(r)] =$$

$$[\mathcal{W}_{m+n}(m+k-1) - \mathcal{W}_m(m+k-1)] \cdot [\mathcal{W}_{m+n}(n+1) - \mathcal{W}_n(n+1)] + \ldots +$$

$$[\mathcal{W}_{m+n}(m+1) - \mathcal{W}_m m+1)] \cdot [\mathcal{W}_{m+n}(n+k-1) - \mathcal{W}_n(n+k-1)] \,.$$

Taking into account $\mathcal{W}_{m+n}(s) = \mathsf{P}(s), \ s \leq m+n$, the above formula can be represented in more convenient form

$$\mathsf{K}_3 = \sum_{k_1, k_2 = 0}^{k_1 + k_2 = k} [\mathsf{P}(m+k_1) - \mathcal{W}_m(m+k_1)] \cdot [\mathsf{P}(n+k_2) - \mathcal{W}_n(n+k_2)] \,.$$

$\underline{n+k \leq r \leq \frac{m+n+k}{2} \leq m+n+k-r \leq m}$

$$\begin{cases} \mathcal{P}_n^m(r) = \mathcal{W}_n(r) \\ \mathcal{W}_{m+n}(r) = \mathcal{W}_m(r) \end{cases}, \begin{cases} \mathcal{P}_n^m(m+n+k-r) = \mathcal{W}_n(m+n+k-r) \\ \mathcal{W}_{m+n}(m+n+k-r) = \mathcal{W}_m(m+n+k-r) \end{cases} \to \mathsf{K}_4 = 0. \ \blacksquare$$

Before enlarging the interval of $s$ where the explicit formulas are valid, we will show that the Proposition 11 is actually more universal than it was formulated in (57).

**Corollary 2**

$\mathcal{D}_n^m(k)$ *does not depend on* $m, n$: $\mathcal{D}_n^m(k) = \mathsf{D}(k) = \mathsf{A}(k) - \mathsf{B}(k)$, $\mathsf{D}(0) = 0$,

$$\mathsf{A}(k) = \sum_{k_1, k_2 = 0}^{k_1 + k_2 = k} \left[ \sum_{q_1 = 0}^{k_1 - 1} \mathsf{P}(q_1) \times \sum_{q_2 = 0}^{k_2 - 1} \mathsf{P}(q_2) \right], \ \mathsf{B}(k) = \sum_{r=0}^{k-1} \mathsf{P}(k-r) \cdot \mathsf{D}(r) \,. \quad (61)$$

*Proof.* First we recall the recursive relation for the restricted partition functions [3]

$$\mathcal{W}_r(s) = \mathcal{W}_{r-1}(s) + \mathcal{W}_r(s-r) \ \to \ \mathcal{W}_{n+k}(s) = \mathcal{W}_n(s) + \sum_{q=0}^{k-1} \mathcal{W}_{n+k-q}(s-n-k+q)$$

and put $s = n + k$. Then we get

$$
\begin{align}
\mathsf{P}(n+k) - \mathcal{W}_n(n+k) &= \sum_{q=0}^{k-1} \mathcal{W}_{n+k-q}(q), \quad \text{and if } n+k \geq 2q \equiv k \leq n+2 \to \\
\mathsf{P}(n+k) - \mathcal{W}_n(n+k) &= \sum_{q=0}^{k-1} \mathsf{P}(q). \tag{62}
\end{align}
$$

The similar formula exists for $k \leq m + 2$. The final steps are clear : since a inequality $k \leq n + 1$ in (57) is more strong than $k \leq n + 2$ in (62) we conclude that the formula (61) must be true for all $k_i \leq k$. ∎

We present here the first ten values of the universal function $\mathsf{D}(k)$:

$$
\begin{align}
&\mathsf{D}(1) = 0, \ \mathsf{D}(2) = 1, \ \mathsf{D}(3) = 3, \ \mathsf{D}(4) = 7, \ \mathsf{D}(5) = 15, \ \mathsf{D}(6) = 25, \\
&\mathsf{D}(7) = 44, \ \mathsf{D}(8) = 75, \ \mathsf{D}(9) = 118, \ \mathsf{D}(10) = 190. \tag{63}
\end{align}
$$

One can encounter a case when one of the restrictions $k \leq n+1$ in (57) is broken while another one is conserved, e.g. $\mathcal{P}_3^{20}(23+6)$. Here we prove

**Proposition 12**

$$
\begin{align}
\mathcal{P}_n^m(m+n+k) &= \mathcal{V}_n^m(m+n+k) + \mathcal{L}_n^m(k), \quad 0 \leq k \leq m-n, \ 2n+1 \geq k > n+1, \\
\mathcal{L}_n^m(k) &= \mathcal{C}_n^m(k) - \mathcal{E}_n^m(k), \ \mathcal{L}_n^m(n+1) = \mathsf{D}(n+1), \\
\mathcal{C}_n^m(k) &= \sum_{r=n}^{k-1} [\mathcal{W}_{m+n}(m+r) - \mathcal{W}_m(m+r)] \cdot [\mathsf{P}(k+n-r) - \mathcal{W}_n(k+n-r)] + \\
&\quad \sum_{\substack{n_1+n_2=n \\ n_1,n_2=0}} [\mathsf{P}(m+n_1) - \mathcal{W}_m(m+n_1)] \cdot [\mathsf{P}(k+n_2) - \mathcal{W}_n(k+n_2)], \\
\mathcal{E}_n^m(k) &= \sum_{r=k-n-1}^{k} \mathsf{P}(r) \cdot \mathsf{D}(k-r) + \sum_{r=1}^{k-n-2} \mathsf{P}(r) \cdot \mathcal{L}_n^m(k-r). \tag{64}
\end{align}
$$

*Proof.* Now we can not use directly the induction based on the decomposition (59) because during this inductive procedure the runing index $r$ passes through the value $n+1$.

The decomposition of Eqn (51) looks like

$$
\mathcal{P}_n^m(m+n+k) = \mathcal{V}_n^m(m+n+k) + \mathsf{K}_5 + \mathsf{K}_6 + \mathsf{K}_7 + \mathsf{K}_8 + \mathsf{K}_9,
$$

$$
\mathsf{K}_5 = \sum_{r=1}^{k-n-2} F_3(r), \ \mathsf{K}_6 = \sum_{r=k-n-1}^{n} F_3(r), \ \mathsf{K}_7 = \sum_{r=n+1}^{k} F_3(r), \ \mathsf{K}_8 = \sum_{r=k+1}^{k+n-1} F_3(r), \ \mathsf{K}_9 = \sum_{r=k+n}^{[\frac{m+n+k}{2}]} F_3(r).
$$

The contributions $\mathsf{K}_i$ in five different regions are
$\underline{1 \leq r \leq k-n-2 < n, \ m+2n+2 \leq m+n+k-r \leq m+n+k-1}$

$$
\begin{cases} \mathcal{W}_m(r) = \mathcal{W}_n(r) = \mathcal{W}_{m+n}(r) = \mathcal{P}_n^m(r) = \mathsf{P}(r) \\ \mathcal{P}_n^m(s) - \mathcal{V}_n^m(s) = \mathcal{L}_n^m(k-r), \ s = m+n+k-r \end{cases} \longrightarrow \mathsf{K}_5 = -\sum_{r=1}^{k-n-2} \mathsf{P}(r) \cdot \mathcal{L}_n^m(k-r).
$$

$\underline{k-n-1 \leq r \leq n, \ m+n+1 < m+k \leq m+n+k-r \leq m+2n+1}$

$$\begin{cases} \mathcal{W}_m(r) = \mathcal{W}_n(r) = \mathcal{W}_{m+n}(r) = \mathcal{P}_n^m(r) = \mathsf{P}(r) \\ \mathcal{P}_n^m(s) - \mathcal{V}_n^m(s) = \mathsf{D}(k-r), \ s = m+n+k-r \end{cases} \longrightarrow \mathsf{K}_6 = -\sum_{k-n-1}^{n} \mathsf{P}(r) \cdot \mathsf{D}(k-r) \ .$$

$\underline{n+1 \leq r \leq k < m, \ m+n \leq m+n+k-r \leq m+k-1 < m+2n+1}$
Here the upper restriction on $k$ appears first time

$$\begin{cases} \mathcal{P}_n^m(r) = \mathcal{W}_n(r), \ \mathcal{W}_{m+n}(r) = \mathcal{W}_m(r) \\ \mathcal{P}_n^m(m+n+k-r) - \mathcal{V}_n^m(m+n+k-r) = \mathsf{D}(k-r) \end{cases} \longrightarrow$$

$$\mathsf{K}_7 = \sum_{r=n}^{k-1}[\mathcal{W}_{m+n}(m+r) - \mathcal{W}_m(m+r)] \cdot [\mathsf{P}(k+n-r) - \mathcal{W}_n(k+n-r)] - \sum_{r=n+1}^{k} \mathsf{P}(r) \cdot \mathsf{D}(k-r) \ .$$

$\underline{k+1 \leq r \leq k+n-1 \leq m-1 < m+1 \leq m+n+k-r \leq m+n-1}$

$$\begin{cases} \mathcal{W}_{m+n}(m+n+k-r) = \mathsf{P}(m+n+k-r), \\ \mathcal{P}_n^m(m+n+k-r) = \mathcal{V}_n^m(m+n+k-r) \end{cases} \text{and} \begin{cases} \mathcal{W}_{m+n}(r) = \mathsf{P}(r), \\ \mathcal{P}_n^m(r) = \mathcal{W}_n^m(r) \end{cases} \longrightarrow$$

In the way, similar to (60), we obtain

$$\mathsf{K}_8 = \sum_{n_1,n_2=0}^{n_1+n_2=n} [\mathsf{P}(m+n_1) - \mathcal{W}_m(m+n_1)] \cdot [\mathsf{P}(k+n_2) - \mathcal{W}_n(k+n_2)] \ .$$

$\underline{k+n \leq r \leq \frac{m+n+k}{2} \leq m+n+k-r \leq m}$

$$\begin{cases} \mathcal{P}_n^m(r) = \mathcal{W}_n(r) \\ \mathcal{W}_{m+n}(r) = \mathcal{W}_m(r) \end{cases}, \begin{cases} \mathcal{P}_n^m(m+n+k-r) = \mathcal{W}_n(m+n+k-r) \\ \mathcal{W}_{m+n}(m+n+k-r) = \mathcal{W}_m(m+n+k-r) \end{cases} \rightarrow \mathsf{K}_9 = 0.$$

We will show that the present Proposition coincides with the Proposition 11 when $k = n+1$. According to the definition (64) we have $\mathsf{K}_5+\mathsf{K}_6 = \sum_1^n F(r)$, $\mathsf{K}_7+\mathsf{K}_8 = \sum_{n+1}^{k+n-1} F(r)$, $\mathsf{K}_9 = 0$ and

$$\begin{cases} \mathcal{C}_n^m(n+1) = \mathsf{A}(n+1) \\ \mathcal{E}_n^m(n+1) = \mathsf{B}(n+1) \end{cases} \rightarrow \mathcal{L}_n^m(n+1) = \mathsf{A}(n+1) - \mathsf{B}(n+1) = \mathsf{D}(n+1) \ . \quad \blacksquare$$

The next restriction ($k \leq m-n$) will be omitted in the following

**Proposition 13**

$$\begin{align}
\mathcal{P}_n^m(m+n+k) &= \mathcal{V}_n^m(m+n+k) + \mathcal{M}_n^m(k), \quad m-n < k \leq n+1, \tag{65} \\
\mathcal{M}_n^m(k) &= \mathcal{K}_n^m(k) - \mathcal{N}_n^m(k), \ \mathcal{M}_n^m(m-n) = \mathsf{D}(m-n), \\
\mathcal{K}_n^m(k) &= \sum_{q=1}^{m-n}[\mathsf{P}(m+k-q) - \mathcal{W}_m(m+k-q)] \cdot [\mathsf{P}(n+q) - \mathcal{W}_n(n+q)] + \\
&\quad \sum_{q=1}^{[\frac{n+k-m}{2}]}[\mathcal{W}_n(m+q)\mathcal{W}_m(n+k-q) - \mathcal{V}_n^m(m+q)\mathsf{P}(n+k-q) + \\
&\quad \mathcal{W}_n(n+k-q)\mathcal{W}_m(m+q) - \mathcal{V}_n^m(n+k-q)\mathsf{P}(m+q)], \\
\mathcal{N}_n^m(k) &= \sum_{r=k-m+n}^{k-2} \mathsf{P}(r) \cdot \mathsf{D}(k-r) + \sum_{r=1}^{k+n-m-1} \mathsf{P}(r) \cdot \mathcal{M}_n^m(k-r).
\end{align}$$

*Proof.* The decomposition of Eqn (51) looks like

$$\mathcal{P}_n^m(m+n+k) = \mathcal{V}_n^m(m+n+k) + \mathsf{K}_{10} + \mathsf{K}_{11} + \mathsf{K}_{12} + \mathsf{K}_{13} + \mathsf{K}_{14},$$

$$\mathsf{K}_{10} = \sum_{1}^{k+n-m-1} F_3(r), \quad \mathsf{K}_{11} = \sum_{r=k+n-m}^{k-2} F_3(r), \quad \mathsf{K}_{12} = \sum_{r=k-1}^{n} F_3(r),$$

$$\mathsf{K}_{13} = \sum_{r=n+1}^{m} F_3(r), \quad \mathsf{K}_{14} = \sum_{r=m+1}^{[\frac{m+n+k}{2}]} F_3(r).$$

$\underline{1 \leq r \leq k+n-m-1 < n, \ 2m+1 \leq m+n+k-r \leq m+n+k-1}$

$$\begin{cases} \mathcal{W}_m(r) = \mathcal{W}_n(r) = \mathcal{W}_{m+n}(r) = \mathcal{P}_n^m(r) = \mathsf{P}(r) \\ \mathcal{P}_n^m(s) - \mathcal{V}_n^m(s) = \mathcal{M}_n^m(k-r), \ s = m+n+k-r \end{cases} \longrightarrow \mathsf{K}_{10} = -\sum_{r=1}^{k+n-m-1} \mathsf{P}(r) \cdot \mathcal{M}_n^m(k-r).$$

$\underline{k+n-m \leq r \leq k-2 < n, \ m+n+2 \leq m+n+k-r \leq 2m < m+n+k < m+2n+1}$

$$\begin{cases} \mathcal{W}_m(r) = \mathcal{W}_n(r) = \mathcal{W}_{m+n}(r) = \mathcal{P}_n^m(r) = \mathsf{P}(r) \\ \mathcal{P}_n^m(s) - \mathcal{V}_n^m(s) = \mathsf{D}(k-r), \ s = m+n+k-r \end{cases} \longrightarrow \mathsf{K}_{11} = -\sum_{r=k-m+n}^{k-2} \mathsf{P}(r) \cdot \mathsf{D}(k-r).$$

$\underline{k-1 \leq r \leq n, \ m+k \leq m+n+k-r \leq m+n+1}$

$$\begin{cases} \mathcal{W}_m(r) = \mathcal{W}_n(r) = \mathcal{W}_{m+n}(r) = \mathcal{P}_n^m(r) = \mathsf{P}(r) \\ \mathcal{P}_n^m(m+n+k-r) = \mathcal{V}_n^m(m+n+k-r) \end{cases} \longrightarrow \mathsf{K}_{12} = 0.$$

$\underline{n+1 \leq r \leq m, \ m < n+k \leq m+n+k-r \leq m+k-1 \leq m+n}$

$$\begin{cases} \mathcal{P}_n^m(r) = \mathcal{W}_n(r) \\ \mathcal{W}_{m+n}(r) = \mathcal{W}_m(r) \end{cases} \text{and} \begin{cases} \mathcal{W}_n(s) - \mathcal{P}_n^m(s) = \mathcal{W}_{m+n}(s) - \mathcal{W}_m(s) \\ m+1 \leq s = m+n+k-r \leq m+n \end{cases} \longrightarrow$$

$$\mathsf{K}_{13} = \sum_{q=1}^{m-n} [\mathsf{P}(m+k-q) - \mathcal{W}_m(m+k-q)] \cdot [\mathsf{P}(n+q) - \mathcal{W}_n(n+q)].$$

$\underline{m+1 \leq r \leq \frac{m+n+k}{2} \leq m+n+k-r \leq n+k-1 < m+n-1}$

$$\begin{cases} \mathcal{P}_n^m(r) = \mathcal{V}_n^m(r) \\ \mathcal{W}_{m+n}(r) = \mathsf{P}(r) \end{cases} \text{and} \begin{cases} \mathcal{P}_n^m(m+n+k-r) = \mathcal{V}_n^m(m+n+k-r) \\ \mathcal{W}_{m+n}(m+n+k-r) = \mathsf{P}(m+n+k-r) \end{cases} \longrightarrow$$

$$\mathsf{K}_{14} = \sum_{q=1}^{[\frac{n+k-m}{2}]} [\mathcal{W}_n(m+q)\mathcal{W}_m(n+k-q) - \mathcal{V}_n^m(m+q)\mathsf{P}(n+k-q) +$$

$$\mathcal{W}_n(n+k-q)\mathcal{W}_m(m+q) - \mathcal{V}_n^m(n+k-q)\mathsf{P}(m+q)].$$

Finally we will show that the present Proposition coincides with the Proposition 11 when $k = m - n$. According to the definition (65) we have $\mathsf{K}_{10} + \mathsf{K}_{11} = \sum_1^{k-2} F(r)$, $\mathsf{K}_{14} = 0$ and

$$\begin{cases} \mathcal{K}_n^m(m-n) = \mathsf{A}(m-n) \\ \mathcal{N}_n^m(m-n) = \mathsf{B}(m-n) \end{cases} \to \mathcal{M}_n^m(m-n) = \mathsf{A}(m-n) - \mathsf{B}(m-n) = \mathsf{D}(m-n). \quad \blacksquare$$

## 4.1 Generalized Gaussian generating function

In this Section we give a proof of (44) for $\mathcal{P}\left[{}^{m_1,m_2}_{n_1,n_2}\right](s)$ and derive its generating function $G(n_1, m_1; n_2, m_2; s)$. First, we reduce a problem to the Diophantine system. For this purpose let us write the counter-partner monomial terms of $\mathsf{S}^{\{m_1,m_2\}}_{S_{n_1} \times S_{n_2}}$ - Eqn (43)

$$T^{s,l}_{\oplus}\left({}^{\{m_1,m_2\}}_{S_{n_1} \times S_{n_2}}\right) = \prod_{r_1=1}^{n_1} I^{\alpha^{s,l}_{r_1}}_{n_1,r_1} \prod_{r_2=1}^{n_2} I^{\beta^{s,l}_{r_2}}_{n_2,r_2} \cdot \lambda^{m_1 n_1 + m_2 n_2 - \sum_{r_1=1}^{n_1} r_1 \alpha^{s,l}_{r_1} - \sum_{r_2=1}^{n_2} r_2 \beta^{s,l}_{r_2}} \quad (66)$$

$$T^{s,l}_{\ominus}\left({}^{\{m_1,m_2\}}_{S_{n_1} \times S_{n_2}}\right) = I^{m_1 - \sum_{r_1=1}^{n_1} \alpha^{s,l}_{r_1}}_{n_1,n_1} I^{m_2 - \sum_{r_2=1}^{n_2} \beta^{s,l}_{r_2}}_{n_2,n_2} \prod_{r_1=1}^{n_1} I^{\alpha^{s,l}_{r_1}}_{n_1,n_1-r_1} \prod_{r_2=1}^{n_2} I^{\beta^{s,l}_{r_2}}_{n_2,n_2-r_2} \cdot \lambda^{\sum_{r_1=1}^{n_1} r_1 \alpha^{s,l}_{r_1} + \sum_{r_2=1}^{n_2} r_2 \beta^{s,l}_{r_2}}$$

or more briefly

$$T^{s,l}_{\oplus}\left({}^{\{m_1,m_2\}}_{S_{n_1} \times S_{n_2}}\right) = T^{s,l}_{\oplus}\left({}^{\{m_1\}}_{S_{n_1}}\right) \times T^{s,l}_{\oplus}\left({}^{\{m_2\}}_{S_{n_2}}\right), \quad T^{s,l}_{\ominus}\left({}^{\{m_1,m_2\}}_{S_{n_1} \times S_{n_2}}\right) = T^{s,l}_{\ominus}\left({}^{\{m_1\}}_{S_{n_1}}\right) \times T^{s,l}_{\ominus}\left({}^{\{m_2\}}_{S_{n_2}}\right),$$

which define the Diophantine system

$$\sum_{r_1=1}^{n_1} r_1 \alpha^{s,l}_{r_1} + \sum_{r_2=1}^{n_2} r_2 \beta^{s,l}_{r_2} = s, \quad \sum_{r_1=1}^{n_1} \alpha^{s,l}_{r_1} \leq m_1, \quad \sum_{r_2=1}^{n_2} \beta^{s,l}_{r_2} \leq m_2. \quad (67)$$

The number $\mathcal{P}\left[{}^{m_1,m_2}_{n_1,n_2}\right](s)$ of positive integer solutions $\{x_{r_1}\}, \{y_{r_2}\}$ of the Diophantine system (67) can be found by successive counting intermediate numbers of solutions

$$\sum_{r_1=1}^{n_1} r_1 \alpha^{s,l}_{r_1} = s - g, \quad \sum_{r_2=1}^{n_2} r_2 \beta^{s,l}_{r_2} = g \longrightarrow \mathcal{P}^{m_1}_{n_1}(s-g) \cdot \mathcal{P}^{m_2}_{n_2}(g), \quad g = 0, 1, \ldots s, \quad (68)$$

that leads to the following generalization of the Cayley-Sylvester Theorem (47)

**Theorem 2**
*Let $\mathcal{P}\left[{}^{m_1,m_2}_{n_1,n_2}\right](s)$ be a number of positive integer solutions $X = \{x_{r_1}\}, Y = \{y_{r_2}\}$ of the system of one Diophantine equation and two Diophantine inequalities*

$$\sum_{r_1=1}^{n_1} r_1 x_{r_1} + \sum_{r_2=1}^{n_2} r_2 y_{r_2} = s, \quad \sum_{r_1=1}^{n_1} x_{r_1} \leq m_1, \quad \sum_{r_2=1}^{n_2} y_{r_2} \leq m_2.$$

*Then $\mathcal{P}\left[{}^{m_1,m_2}_{n_1,n_2}\right](s)$ is generated by the Gaussian polynomial $G(n_1, m_1; n_2, m_2; s)$ of the finite order $n_1 m_1 + n_2 m_2$*

$$G(n_1, m_1; n_2, m_2; t) = G(n_1, m_1; t) \cdot G(n_2, m_2; t) = \sum_{s \geq 0}^{n_1 m_1 + n_2 m_2} \mathcal{P}\left[{}^{m_1,m_2}_{n_1,n_2}\right](s) \cdot t^s, \quad (69)$$

$$\mathcal{P}\left[{}^{m_1,m_2}_{n_1,n_2}\right](s) = \sum_{s_1,s_2=0}^{s} \mathcal{P}^{m_1}_{n_1}(s_1) \cdot \mathcal{P}^{m_2}_{n_2}(s_2), \quad s_1 + s_2 = s,$$

*where partial Gaussian polynomial $G(n_j, m_j; t)$ is defined in (48).*

*Proof.* The proof follows immediately after summation in (68) over $g$ and comparison with serial expansion of $G(n_1, m_1; t) \cdot G(n_2, m_2; t)$ in accordance with (48). ∎

We finish this Section with a simple extension of the Theorem 2 onto system of one Diophantine equation and $k$ Diophantine inequalities. This generalization concerned with skew-reciprocal Eqn which is built upon the basis of polynomial invariants of the direct product of finite number of symmetric groups $G = S_{n_1} \times \ldots \times S_{n_k}$.

**Corollary 3**
Let $\mathcal{P}\begin{bmatrix} m_1,...,m_k \\ n_1,...,n_k \end{bmatrix}(s)$ be the number of positive integer solutions $X_j = \{x_{j,r_j}\}$ of the system of one Diophantine equation and $k$ Diophantine inequalities

$$\sum_{r_1=1}^{n_1} r_1\, x_{1,r_1} + \ldots + \sum_{r_k=1}^{n_k} r_k\, x_{k,r_k} = s\,, \quad \sum_{r_j=1}^{n_j} x_{j,r_j} \leq m_j\,, \quad j=1,...,k\,. \qquad (70)$$

Then $\mathcal{P}\begin{bmatrix} m_1,...,m_k \\ n_1,...,n_k \end{bmatrix}(s)$ is generated by the Gaussian polynomial $G(n_1, m_1; ...; n_k, m_k; t)$ of finite order $NM = \sum_{j=1}^{k} n_j m_j$

$$G(n_1, m_1;\, \ldots\, ; n_k, m_k; t) = \prod_{j=1}^{k} G(n_j, m_j; t) = \sum_{s \geq 0}^{NM} \mathcal{P}\begin{bmatrix} m_1,...,m_k \\ n_1,...,n_k \end{bmatrix}(s) \cdot t^s\,, \qquad (71)$$

$$\mathcal{P}\begin{bmatrix} m_1,...,m_k \\ n_1,...,n_k \end{bmatrix}(s) = \sum_{s_1,...,s_k=0}^{s} \prod_{j=1}^{k} \mathcal{P}_{n_j}^{m_j}(s_j)\,, \quad s_1 + ... + s_k = s\,.$$

The proof is similar to the previous in Theorem 2.

## 5 Symmetric unimodality index

In this Section we will find the unimodality indices $\mu\left\{\mathsf{R}_{S_n}^{\{m\}}\right\}$ and $\mu\left\{\mathsf{S}_{S_n}^{\{m\}}\right\}$ based on the definition (11). First, we establish a fact that will be useful, namely

$$\mu\left\{\mathsf{T1}_{S_n}^{\{m\}}\right\} + \mu\left\{\mathsf{T2}_{S_n}^{\{m\}} + \mathsf{T3}_{S_n}^{\{m\}}\right\} + \mu\left\{\mathsf{T4}_{S_n}^{\{m\}}\right\} = \sum_{s \geq 0}^{nm} \mathcal{P}_n^m(s) = \binom{m+n}{n}\,. \qquad (72)$$

Indeed, the first equation follows from the definition (11) of unimodality indices while the second one comes from an accurate calculation of the limit in (48)

$$\sum_{s \geq 0}^{nm} \mathcal{P}_n^m(s) = \lim_{t \to 1} G(n, m; t) = \frac{(m+n)!}{m!\, n!}\,.$$

By the other hand unimodality indices are related to the conformal partitions $\mathcal{P}_n^m(s)$ in the following way

$$\mu\left\{\mathsf{R}_{S_n}^{\{m\}}\right\} = \sum_{s=1}^{[\frac{mn}{2}]} \mathcal{P}_n^m(s)\,, \quad \mu\left\{\mathsf{R}_{S_n}^{\{m\}}\right\} - \mu\left\{\mathsf{S}_{S_n}^{\{m\}}\right\} = \begin{cases} 0\,, & \text{if } mn = 2k+1\,, \\ \mathcal{Q}(m,n)\,, & \text{if } mn = 2k\,, \end{cases} \qquad (73)$$

where $\mathcal{Q}(m,n)$ is a number of counter-partner terms which annihilate in $\mathsf{T2}_{S_n}^{\{m\}}(\lambda_n, x_i) - \mathsf{T3}_{S_n}^{\{m\}}(\lambda_n, x_i)$. E.g.

$$\mathsf{T2}_{S_3}^{\{4\}}(\lambda_n, x_i) + \mathsf{T3}_{S_3}^{\{4\}}(\lambda_n, x_i) = c_{6,1}(I_{3,1}^3 I_{3,3} + I_{3,2}^3) + c_{6,2} I_{3,1}^2 I_{3,2}^2 + c_{6,3} I_{3,1} I_{3,2} I_{3,3} + c_{6,4} I_{3,3}^2\,,$$
$$\mathsf{T2}_{S_3}^{\{4\}}(\lambda_n, x_i) - \mathsf{T3}_{S_3}^{\{4\}}(\lambda_n, x_i) = C_{6,1}(I_{3,1}^3 I_{3,3} - I_{3,2}^3)\,.$$

Eqns (72) and (73) give :

$mn = 2f+1$

$$\mu\left\{\mathsf{R}_{S_n}^{\{m\}}\right\} = \mu\left\{\mathsf{S}_{S_n}^{\{m\}}\right\} = \frac{1}{2}\binom{m+n}{n}, \qquad (74)$$

$mn = 2f$

$$\mu\left\{\mathsf{R}_{S_n}^{\{m\}}\right\} = \frac{1}{2}\binom{m+n}{n} + \frac{1}{2}\mathcal{Q}(m,n)\ , \quad \mu\left\{\mathsf{S}_{S_n}^{\{m\}}\right\} = \frac{1}{2}\binom{m+n}{n} - \frac{1}{2}\mathcal{Q}(m,n)\ . \qquad (75)$$

We will show now that $\mathcal{Q}(m,n)$ is related to the *Catalan partitions* if $mn$ is an even integer. Really, let $mn = 2f$ and the annihilation condition for the counter-partner terms at $s^* = \frac{mn}{2}$ in (9) yields

$$\sum_{r=1}^{n} \alpha_r^{s^*,l} \ln\left(\frac{I_{n,r}}{I_{n,n-r}}\right) = \left(m - \sum_{r=1}^{n} \alpha_r^{s^*,l}\right)\ln I_{n,n}\ , \ \Rightarrow\ \alpha_r^{s^*,l} = \alpha_{n-r}^{s^*,l}\ ,\ m = 2\,\alpha_n^{s^*,l} + \sum_{r=1}^{n-1}\alpha_r^{s^*,l}$$

due to algebraic independence of basic invariants. The latter Diophantine Eqn will be considered in both cases

$$1)\ n = 2u+1\ ,\ \ \frac{m}{2} = \alpha_{2u+1}^{s^*,l} + \sum_{r=1}^{u}\alpha_r^{s^*,l}\ ;\ \ 2)\ n = 2u\ ,\ \ m = 2\left(\alpha_{2u}^{s^*,l} + \sum_{r=1}^{u-1}\alpha_r^{s^*,l}\right) + \alpha_u^{s^*,l}.\ (76)$$

We recall a simple formula for Catalan partitions [14] :

the Diophantine equation $x_1 + x_2 + \ldots + x_n = m$ has $\binom{m+n-1}{m}$ sets of non-negative solutions, the repetitions are not excluded.

Now the both cases (76) are resolved easily

$$\left.\begin{array}{lll} 1)\ \ n = 2u+1\ ,\ & m = 2v\ ,\ & v = \alpha_{2u+1}^{s^*,l} + \sum_{r=1}^{u}\alpha_r^{s^*,l} \\ 2a)\ n = 2u\ ,\ & m = 2v\ ,\ & v = \alpha_{2u}^{s^*,l} + \sum_{r=1}^{u-1}\alpha_r^{s^*,l} + \frac{\alpha_u^{s^*,l}}{2} \\ 2b)\ n = 2u\ ,\ & m = 2v+1\ ,\ & v = \alpha_{2u}^{s^*,l} + \sum_{r=1}^{u-1}\alpha_r^{s^*,l} + \frac{\alpha_u^{s^*,l}-1}{2} \end{array}\right\} \rightarrow \mathcal{Q}(m,n) = \binom{v+u}{v}\ (77)$$

where the fractions in (77) are non-negative integer numbers. The case (74) also can be considered as a trivial solution of Diophantine Eqn. Indeed,

$$n = 2u+1\ ,\ m = 2v+1\ \rightarrow\ 2v+1 = 2\alpha_{2u+1}^{s^*,l} + 2\sum_{r=1}^{u}\alpha_r^{s^*,l}\ \rightarrow\ \mathcal{Q}(m,n) = 0\ ,\qquad (78)$$

since the Diophantine Eqn in (78) has no solution. The formulas (75), (77) present the unimodality indices.

In the same manner we can find the unimodality indices for $\mathsf{S}\begin{bmatrix}m_1,\ldots,m_k\\S_1\times\ldots\times S_k\end{bmatrix}$- and $\mathsf{R}\begin{bmatrix}m_1,\ldots,m_k\\S_1\times\ldots\times S_k\end{bmatrix}$- polynomials

$$\mu\left\{\mathsf{R}\begin{bmatrix}m_1,\ldots,m_k\\S_1\times\ldots\times S_k\end{bmatrix}\right\} = \frac{1}{2}\prod_{j=1}^{k}\binom{m_j+n_j}{n_j} + \frac{1}{2}\mathcal{Q}(n_1,m_1;\ldots;n_k,m_k)\ ,$$

$$\mu\left\{\mathsf{S}\begin{bmatrix}m_1,\ldots,m_k\\S_1\times\ldots\times S_k\end{bmatrix}\right\} = \frac{1}{2}\prod_{j=1}^{k}\binom{m_j+n_j}{n_j} - \frac{1}{2}\mathcal{Q}(n_1,m_1;\ldots;n_k,m_k)\ ,\qquad (79)$$

where the product in (79) arises due to calculation of the limit of $G(n_1, m_1; ...; n_k, m_k; t)$ when $t \to 1$ and $\mathcal{Q}(n_1, m_1; ...; n_k, m_k)$ is a number of counter-partner terms

$$T_\oplus^{s,l}\left(\begin{smallmatrix}\{m_1,...,m_k\}\\S_{n_1}\times...\times S_{n_k}\end{smallmatrix}\right) = \prod_{j=1}^{k} T_\oplus^{s,l}\left(\begin{smallmatrix}\{m_j\}\\S_{n_j}\end{smallmatrix}\right) = \prod_{j=1}^{k}\prod_{r_j=1}^{n_j} I_{n_j,r_j}^{\alpha_{r_j}^{s,l}} \cdot \lambda^{m_j n_j - \sum_{r_j=1}^{n_j} r_j \alpha_{r_j}^{s,l}},$$

$$T_\ominus^{s,l}\left(\begin{smallmatrix}\{m_1,...,m_k\}\\S_{n_1}\times...\times S_{n_k}\end{smallmatrix}\right) = \prod_{j=1}^{k} T_\ominus^{s,l}\left(\begin{smallmatrix}\{m_j\}\\S_{n_j}\end{smallmatrix}\right) = \prod_{j=1}^{k} I_{n_j,n_j}^{m_j - \sum_{r_j=1}^{n_j} \alpha_{r_j}^{s,l}} \prod_{r_j=1}^{n_j} I_{n_j,n_j-r_j}^{\alpha_{r_j}^{s,l}} \cdot \lambda^{\sum_{r_j=1}^{n_j} r_j \alpha_{r_j}^{s,l}},$$

which annihilate. This gives the following implication at $s^\bullet = \frac{1}{2}\sum_{j=1}^{k} n_j, m_j$:

$$\alpha_{r_j}^{s^\bullet,l} = \alpha_{n_j-r_j}^{s^\bullet,l}, \quad m_j = 2\,\alpha_{n_j}^{s^\bullet,l} + \sum_{r_j=1}^{n_j-1} \alpha_{r_j}^{s^\bullet,l} \longrightarrow \mathcal{Q}(n_1, m_1; ...; n_k, m_k) = \prod_{j=1}^{k} \mathcal{Q}(n_j, m_j), \quad (80)$$

where $\mathcal{Q}(n_j, m_j)$ are defined in (77, 78). In order to make (80) consistent with (75) we must supplement the definition of $\mathcal{Q}(n, m)$ with following : $\mathcal{Q}(n, 0) = \mathcal{Q}(0, m) = 1$.

Taking in mind (77, 78, 80) we get

$$\mathcal{Q}(n_1, m_1; ...; n_k, m_k) = 0, \quad \text{if there exists at least one pair } n_j, m_j \text{ that } n_j m_j = 2f_j + 1,$$

$$\mathcal{Q}(n_1, m_1; ...; n_k, m_k) = \prod_{j=1}^{k} \binom{v_j + u_j}{v_j}, \quad \text{if all pairs } n_j, m_j \text{ satisfy } n_j m_j = 2f_j, \quad (81)$$

and $n_j, m_j$ could be represented in one of three different ways (77).

# 6  Conclusion

1. We have considered the self-dual symmetric polynomials — reciprocal $\mathsf{R}_{S_n}^{\{m\}}(\lambda_n, x_i)$ and skew-reciprocal $\mathsf{S}_{S_n}^{\{m\}}(\lambda_n, x_i)$ based on the polynomial invariants of the symmetric group $S_n$ and have studied their properties. They form an infinite commutative semigroup. Real solutions $\lambda_n(x_i)$ of corresponding algebraic Eqns are homogeneous of 1-st order algebraic functions, which are invariant both upon the action of the symmetric group $S_n$, permuting $n$ non-negative variables $x_i$, and conformal group $\mathsf{W}$, inverting both function $\lambda_n$ and the variables $x_i$. We have also discussed some kinds of skew-reciprocal Eqns when one can arrive at their exact solutions $\lambda_n(x_i)$ or at least establish the upper and lower bounds for $\lambda_n(x_i)$.

2. We have discussed the existence of other finite groups $G$ distinguished from symmetric $S_n$ that makes it possible to build out the reciprocal $\mathsf{R}_G$ – and skew-reciprocal $\mathsf{S}_G$ – Eqns on the basis of homogeneous polynomial invariants $I_{d_k}(G)$ of such group. The necessary condition of their existence imposes the duality relation onto degrees $d_k$ of invariants $I_{d_k}(G)$, which compose those Eqns. This condition requires also a monomiality of the invariant $I_{d_{max}}(G)$ with highest degree $d_{max}$. We have shown that some classical groups — cyclic $\mathcal{Z}_n$, alternating $\mathcal{A}_n$, Coxeter $\mathsf{B_n}$, $\mathsf{D_n}$ and symmetric groups' product $S_{n_1} \times ... \times S_{n_k}$ — give rise to $\mathsf{R}_G$ – and $\mathsf{S}_G$ – Eqns while some other groups —- Coxeter $\mathsf{I_{2,m}}$, $m \ne 2, 4$, $\mathsf{G_2}$ and $\mathsf{H_3}$ do not. The similar question concerned with other classical groups remains open.

3. We have derived the analytic expression (52) for conformal partition function $\mathcal{P}_n^m(s)$ based on its relation to usual restricted partition function $\mathcal{W}_n(s)$. The former function is concerned with a special sort of algebraic polynomials — reciprocal $\mathsf{R}_{S_n}^{\{m\}}$ and skew-reciprocal $\mathsf{S}_{S_n}^{\{m\}}$ based on the polynomial invariants of symmetric group $S_n$. We have given also the generalization of Cayley-Sylvester Theorem for a system of one Diophantine equation and $k$ Diophantine inequalities and have derived the corresponding Gaussian polynomial $G(n_1, m_1; ...; n_k, m_k; t)$ and conformal partition function $\mathcal{P}\begin{bmatrix} m_1,...,m_k \\ n_1,...,n_k \end{bmatrix}(s)$.

4. We have found simple expressions (75), (77) for unimodality indices $\mu\left\{\mathsf{R}_{S_n}^{\{m\}}\right\}$, $\mu\left\{\mathsf{S}_{S_n}^{\{m\}}\right\}$ of these polynomials. The unimodality indices $\mu\left\{\mathsf{R}\begin{bmatrix} m_1,...,m_k \\ S_1\times...\times S_k \end{bmatrix}\right\}$ and $\mu\left\{\mathsf{S}\begin{bmatrix} m_1,...,m_k \\ S_1\times...\times S_k \end{bmatrix}\right\}$ of generalized reciprocal and skew-reciprocal polynomials were found also.

# 7  Acknowledgement

I am thankful to Yu. Lyubich for helpful discussion. The encouraging comments of Prof. G. E. Andrews are highly appreciated. This work was began during my stay at the Isaac Newton Institute for Mathematical Sciences which hospitality are highly appreciated.

The research was supported in part by grants from the Tel Aviv University Research Authority and the Gileadi Fellowship program of the Ministry of Absorption of the State of Israel.

# A  Triangular Toeplitz inversion formula

We will study a linear convolution equations with a triangular Toeplitz matrix

$$\mathcal{P}(g) = T(g) + \sum_{s=0}^{g-1} \mathcal{P}(s) \cdot U(g-s) \,, \quad \mathcal{P}(0) = T(0) \,, \tag{A1}$$

where two known functions $T(g), U(g)$ and unknown one $\mathcal{P}(g)$ are considered only on the non-negative integers. The successive recursion of (A1) gives

$$\begin{aligned}\mathcal{P}(g) &= T(g) + T(g-1) \cdot U(1) + \sum_{s=0}^{g-2} \mathcal{P}(s) \cdot \{\, U(g-s) + U(1) \cdot U(g-1-s) \,\} = \\ &= T(g) + T(g-1) \cdot U(1) + T(g-2) \cdot [\, U(2) + U^2(1) \,] + \\ &\quad \sum_{s=0}^{g-3} \mathcal{P}(s) \cdot \{U(g-s) + U(1) \cdot U(g-1-s) + [\, U(2) + U^2(1) \,] \cdot U(g-2-s)\} \,.\end{aligned}$$

Now one can arrive by induction at

$$\mathcal{P}(g) = \sum_{r=0}^{k-1} T(g-r) \cdot \Phi_r(U) + \sum_{s=0}^{g-k} \mathcal{P}(s) \cdot \sum_{r=0}^{k-1} U(g-r-s) \cdot \Phi_r(U) \,, \quad 1 \le k \le g \,, \tag{A2}$$

where polynomials $\Phi_r(U)$ are related to the restricted partition number $\mathcal{W}_g(r)$ of positive integer $r$ into parts, none of which exceeds $g$. They are built out of $\mathcal{W}_r(r)$ algebraically independent terms with coefficients the sum of which is equal $2^{r-1}$

$$\Phi_r(U) = \sum_{q=1}^{r} q! \sum_{\{q_l\}}^{q=const} \frac{1}{q_1! \cdot q_2! \cdot \ldots \cdot q_r!} \prod_{l=1}^{r} U^{q_l}(l) \,, \quad \sum_{l=1}^{r} q_l \cdot l = r \,, \quad \sum_{l=1}^{r} q_l = q \,,$$

i.e.

$$\begin{aligned}\Phi_0(U) &= 1 \,, \quad \Phi_1(U) = U(1) \,, \quad \Phi_2(U) = U(2) + U^2(1) \,, \\ \Phi_3(U) &= U(3) + 2\, U(2) \cdot U(1) + U^3(1) \,, \\ \Phi_4(U) &= U(4) + 2\, U(3) \cdot U(1) + U^2(2) + 3\, U(2) \cdot U^2(1) + U^4(1) \,, \\ \Phi_5(U) &= U(5) + 2\, U(4) \cdot U(1) + 2\, U(3) \cdot U(2) + 3\, U(3) \cdot U^2(1) + 3\, U^2(2) \cdot U(1) + \\ &\quad 4\, U(2) \cdot U^3(1) + U^5(1) \,, \\ \Phi_6(U) &= U(6) + 2\, U(5) \cdot U(1) + 2\, U(4) \cdot U(2) + U^2(3) + 3\, U(4)\, U^2(1) + U^3(2) + \\ &\quad 6\, U(3)\, U(2)\, U(1) + 4\, U(3)\, U^3(1) + 6\, U^2(2)\, U^2(1) + 5\, U(2)\, U^4(1) + U^6(1) \,.\end{aligned} \tag{A3}$$

After putting $k = g$ in (A2) we will get finally

$$\mathcal{P}(g) = \sum_{r=0}^{g-1} [\, T(g-r) + U(g-r) \,] \cdot \Phi_r(U) \,. \tag{A4}$$

Let us compare now Eqn (A1) with Eqn (51) :

$$T(g) = \sum_{s=0}^{g} \mathcal{W}_n(s) \cdot \mathcal{W}_m(g-s) \,, \quad U(s) = -\mathcal{W}_{n+m}(s) \,.$$

Then we immediately arrive at (52). One can verify that the formulas (54, 57) can be obtained using (52).